\documentclass[a4paper]{article}
\usepackage{amssymb}
\usepackage{amsthm}
\usepackage{amsmath}
\usepackage{latexsym}
\usepackage[T1]{fontenc}
\usepackage[latin1]{inputenc}

\def \eop {\hbox{}\nobreak\hfill \vrule width 2.0mm height 1.8mm depth 0mm
\par \goodbreak \smallskip}

\newcommand{\integ}[2]{\displaystyle \int_{#1}^{#2}}
\newcommand{\dint}{\displaystyle\int}

\begin{document}
\newtheorem{definition}{Definition}[section]
\newtheorem{theorem}{Theorem}[section]
\newtheorem{proposition}{Proposition}[section]
\newtheorem{lemma}{Lemma}[section]
\newtheorem{remark}{Remark}[section]
\newtheorem{corollary}{Corollary}[section]
\newtheorem{example}{Example}[section]
\def \ce{\centering}
\def \bop {\noindent\textbf{Proof}}
\def \eop {\hbox{}\nobreak\hfill
\vrule width 2mm height 2mm depth 0mm
\par \goodbreak \smallskip}
\def \R{I\!\!R}
\def \N{I\!\!N}
\def \P{\mathbb{P}}
\def \E{I\!\!E}
\def \T{\mathbb{T}}
\def \H{\mathbb{H}}
\def \L{\mathbb{L}}
\def \Z{\mathbb{Z}}
\def \bf{\textbf}
\def \it{\textit}
\def \sc{\textsc}
\def \ni {\noindent}
\def \sni {\ss\ni}
\def \bni {\bigskip\ni}
\def \ss {\smallskip}
\def \F{\mathcal{F}}
\def \g{\mathcal{g}}
\def \eop {\hbox{}\nobreak\hfill
\vrule width 2mm height 2mm depth 0mm
\par \goodbreak \smallskip}

\title{ Generalized BSDE with 2-Reflecting Barriers and Stochastic Quadratic Growth.}
\author{E. H. Essaky \quad \quad M. Hassani\\\\
Universit\'{e} Cadi Ayyad\\ Facult\'{e} Poly-disciplinaire\\
Laboratoire de Modélisation et Combinatoire\\
D\'{e}partement de Math\'{e}matiques et d'Informatique\\ B.P 4162, Safi, Maroc.\\
e-mails: essaky@ucam.ac.ma \hspace{.2cm}medhassani@ucam.ac.ma}
\date{}
\maketitle \footnotetext[1]{This work is supported by Hassan II
Academy of Science and technology and Action Int\'egr\'ee
MA/10/224.}


\begin{abstract}
We study the existence of a solution for a one-dimensional
generalized backward stochastic differential equation with two
reflecting barriers (GRBSDE for short) under assumptions on the
input data which are weaker than that on the current literature. In
particular, we construct a maximal solution for such a GRBSDE when
the terminal condition $\xi$ is only ${\cal F}_T-$measurable and the
driver $f$ is continuous with general growth with respect to the
 variable $y$ and stochastic quadratic growth with respect to the
 variable $z$ without assuming any $P-$integrability conditions.

 The work is suggested by
the interest the results might have in Dynkin game problem and
American game option.
\end{abstract}

\ni \textbf{Keys Words:} Reflected backward stochastic differential
equation; stochastic quadratic growth; comparison theorem;
exponential transformation.
\medskip

\ni \textbf{AMS Classification}\textit{(1991)}\textbf{: }60H10,
60H20.
\section{Introduction}

\indent Backward stochastic differential equations (BSDEs for short)
have been introduced long time ago by J. B. Bismut \cite{bis} both
as the equations for the adjoint process in the stochastic version
of Pontryagin maximum principle as well as the model behind the
Black and Scholes formula for the pricing and hedging of options in
mathematical finance. However the first published paper on nonlinear
BSDEs appeared only in 1990, by Pardoux and Peng \cite{pp1}. A
solution for such an equation is a couple of adapted processes $(Y,
Z)$ with values in $\R\times\R^d$ satisfying
\begin{equation}\label{equa0}
 Y_t = \xi + \dint_t^T f(s,Y_s,Z_s)ds - \dint_t^T Z_s dB_s,
\quad\quad 0\leq t\leq T.
\end{equation}
In \cite{pp1}, the authors have proved the existence and uniqueness
of the solution under conditions including basically the Lipschitz
continuity of the generator $f$.

Later on, the study of BSDEs has been motivated by their many
applications in mathematical finance, stochastic control and the
second order PDE theory (see, for example, \cite{KPQ, HL, PP2, pp1,
BH, BHP, Kob} and the references therein).


The notion of BSDE with two reflecting barriers has been first
introduced by Civitanic and Karatzsas \cite{CK}. A solution for such
an equation, associated with a coefficient $f$; terminal value $\xi$
and two barriers $L$ and $U$, is a quadruple of processes $(Y, Z,
K^+, K^-)$ with values in $\R\times\R^d\times \R_+\times \R_+$
satisfying:
\begin{equation}
\label{eq000} \left\{
\begin{array}{ll}
(i) &
 Y_{t}=\xi
+\integ{t}{T}f(s,Y_{s},Z_{s})ds +\integ{t}{T}dK_{s}^+
-\integ{t}{T}dK_{s}^- -\integ{t}{T}Z_{s}dB_{s}\,, t\leq T,
\\ (ii)&
\forall  t\leq T,\,\, L_t \leq Y_{t}\leq U_{t},\quad \\
(iii)& \integ{0}{T}( Y_{t}-L_{t}) dK_{t}^+= \integ{0}{T}(
U_{t}-Y_{t}) dK_{t}^-=0,\,\, \mbox{a.s.}, \\ (iv) &  K_0^+ =K_0^-
=0, \,\,\,\, K^+, K^-, \,\,\mbox{are continuous nondecreasing}.
\end{array}
\right.
 \end{equation}
 Here the solution process $Y$ has to remain between $L$ and $U$ due to the
 cumulative actions of processes $K^+$ and $K^-$. In the case
of a uniformly Lipschitz coefficient $f$ and a square terminal
 condition $\xi$ the existence and uniqueness of a solution have been proved
 when the barriers $L$ and $U$ are either regular or satisfy Mokobodski's condition which, roughly speaking,
 turns out into the existence of a difference of nonnegative
 supermartingales between $L$ and $U$. It has been shown also in \cite{CK} that the solution coincides with the
value of a stochastic Dynkin game of optimal stopping. The link
between obstacle PDEs and RBSDEs has been given in Hamad\`ene and
Hassani \cite{HH}. In Hamad\`ene \cite{Ham1}, applications of RBSDEs
to Dynkin games theory as well as to American game option are given.
We should mention here that, in the case of one-side reflected BSDE,
Kobylanski, Lepeltier, Quenez and Torres \cite{KLQT} provide
existence of solution when the coefficient $f$ is continuous, has a
superlinear growth in $y$ and quadratic growth in $z$. They also
give a characterization of the solution as the value function of an
optimal stopping time problem.

A natural question is then arises : are there any  weaker conditions
under which the RBSDE (\ref{eq000}) has a solution? This question
has attracted many authors in this area. For example, when the
generator $f$ is only continuous there exists a solution to BSDE
(\ref{eq000})under one of the following group of conditions :\\

$\bullet$ $\xi$ is square integrable, $f$ has a uniform linear
growth in $y$ and $z$, i.e. there exists a constant $C$ such that
$|f(t,\omega, y, z)|\leq C(1+|y|+ |z|)$, and one of the barriers has
to be regular, \it{e.g.} has to be
semi-martingale (see Hamad\`ene et \it{al} \cite{HLM}).\\

$\bullet$ $\xi$ is bounded, $f$ has a general growth in $y$ and
quadratic growth in $z$, i.e. there exist a constant $C$ and
positive function $\phi$ which is bounded on compacts such that
$|f(t,\omega, y, z)|\leq C(1+\phi(|y|)+ |z|^2)$, and the barriers
satisfy the Mokobodski's condition (see Bahlali et \it{al}
\cite{BHM}).\\

$\bullet$ $\xi$ is square integrable, $f$ has a uniform linear
growth in $y$ and $z$ and the barriers are square integrable and
completely separated \it{i.e.} $L_t < U_t,\,\, \forall t\in [0,T]$
(see Hamad\`ene and Hassani \cite{HH}).\\

$\bullet$ $\xi$ is square integrable, $f$ linearly increasing and
the barriers are such that $L < U$ on $[0,T)$ and there exists a
continuous semimartingale between $L$ and $U$ (see Lepeltier and San
Martin \cite{LS2}).


%
\vspace{0.5cm} The main objective of this work is to extend and
improve the existence conditions of a solution for GRBSDE
(\ref{eq000}). So the new features here are. Firstly, the generator
$f$ is continuous with general growth with respect to the variable
$y$ and stochastic quadratic growth with respect to the variable $z$
of the form $C_s(\omega)\mid z\mid^2$ instead of $C\mid z\mid^2$ as
usually done. Secondly, instead of assuming the Mokobodski's
condition on the barriers $L$ and $U$, we suppose only that there
exists a semimartingale between them. Thirdly, we do not assume any
assumptions on the $P$-integrability on the input data. We present
also in the appendix a
comparison theorem under general assumptions on the coefficients.\\

By means of an exponential change, the proof of our main result
consists in establishing first a correspondence between our GRBSDE
and another GRBSDE whose coefficients are more tractable. We show
that the existence of solutions for our initial GRBSDE is equivalent
to the existence of solutions for the auxiliary GRBSDE. Since the
integrability conditions on parameters are weaker, we make use of
approximations and truncations to establish the existence result for
the auxiliary GRBSDE. The final step consists in justifying the
passage to the limit and in identifying the limit as the solution of
the auxiliary GRBSDE.



 Let us describe our plan. First, most of the material used
in this paper is defined in Section 2, an exponential transformation
for our GBSDE with two reflecting barriers is also given. In Section
3, with the help of the comparison theorem and using an
approximation technique, we prove the existence of a maximal
solution for the transformed BSDE and then equivalently the
existence of maximal solution for our GBSDE with two reflecting
barriers.  Finally, in appendix, we give a comparison theorem for a
general GBSDE with two reflecting barriers as well as the existence
and uniqueness of a solution of Equation (\ref{eq000}) when the
coefficients $f$ and $g$ are bounded Lipschitz functions and the
processes $L, U, \xi, A$ and $R$ are bounded.

\section{Problem formulation, assumptions and exponential transformation for GBSDE}
In this section, we collect some preliminary results which will be
useful in the sequel.
\subsection{Assumptions and remarks}
Let $(\Omega, {\cal F}, ({\cal F}_t)_{t\leq T}, P)$ be a stochastic
basis on which is defined a $d-$dimensional Brownian motion
$(B_t)_{t\leq T}$ such that $({\cal F}_t)_{t\leq T}$ is the natural
filtration of $(B_t)_{t\leq T}$ and ${\cal F}_0$ contains all
$P$-null sets of $\cal F$.  Note that $({\cal F}_t)_{t\leq T}$
satisfies the usual conditions, \it{i.e.} it is right continuous and
complete.
\medskip\\ Let us now introduce the following notations :
\medskip

$\bullet$ $\cal P$ the sigma algebra of ${\cal F}_t$-progressively
measurable sets on $\Omega\times[0,T].$

$\bullet$ ${\cal C}$ the set of $\R$-valued $\cal P$-measurable
continuous processes $(Y_t)_{t\leq T}$.

$\bullet$ ${\cal L}^{2,d}$ the set of $\R^d$-valued and $\cal
P$-measurable processes $(Z_t)_{t\leq T}$ such that
$$\integ{0}{T}|Z_s|^2ds<\infty, P- a.s.$$

$\bullet$ ${\cal M}^{2,d}$ the set of $\R^d$-valued and $\cal
P$-measurable processes $(Z_t)_{t\leq T}$ such that
$$\E\integ{0}{T}|Z_s|^2ds<\infty.$$

$\bullet$ ${\cal K}$ the set of ${\cal P}$-measurable continuous
nondecreasing
processes $(K_t)_{t\leq T}$ such that $K_0 = 0$ and $K_T <+\infty,$ $P$-- a.s.\\

$\bullet$ ${\cal K}-{\cal K}$ the set of ${\cal P}$-measurable
and continuous processes $(V_t)_{t\leq T}$ such that there exist $V^+, V^-\in{\cal K}$ satisfying : $V = V^+ -V^-$.\\

$\bullet$ ${\cal K}^2$ the set of ${\cal P}$-measurable continuous
nondecreasing processes $(K_t)_{t\leq T}$ such that $K_0 = 0$ and
$\E K_T^2 <+\infty$.

$\bullet$ $A$ and $R$ are two processes in $\cal K$ and $\cal K
-\cal K$ respectively.\\

\ni The following notations are also needed :\\

$\bullet$ For a set $B$, we denote by $B^c$ the complement of $B$
and $1_B$ denotes the indicator of $B$.

$\bullet$ For each $(a, b)\in\R^2$, $a\wedge b = \min(a, b)$\,\, and
\,\,$a\vee b = \max(a, b)$.

$\bullet$ For all $(a, b, c)\in\R^3$ such that $a\leq c$,\,\, $a\vee
b\wedge c = \min(\max(a,b), c) = \max(a, \min(c,b))$.

\medskip
\ni To give conditions under which solutions to a GRBSDE exist, we
should first give the following definition.
\begin{definition}
Let $K^1$ and $K^2$ be two processes in ${\cal K}$. We say that :
\begin{enumerate}
\item $K^1$ and $K^2$ are singular if and only if there exists a set
$D\in {\cal P}$ such that
$$
\E\dint_0^T 1_D(s,\omega) dK^1_s(\omega) = \E\dint_0^T
1_{D^c}(s,\omega) dK^2_s(\omega) =0.
$$
 This is denoted by $dK^1 \perp dK^2$.
\item $dK^1 \leq dK^2$ if and only if for each set $B\in {\cal P}$
$$
\E\dint_0^T 1_B(s,\omega) dK^1_s(\omega)\leq \E\dint_0^T
1_{B}(s,\omega) dK^2_s(\omega), \quad \it{i.e.}\quad K_t^1
-K_s^1\leq K_t^2-K_s^2,\,\,\, \forall s\leq t\quad P-a.s.
$$
In this case $\dfrac{dK^1}{dK^2}$ denotes a ${\cal P}-$measurable
Radon-Nikodym density of $dK^1$ with respect to $dK^2$ which
satisfies
$$
0\leq \dfrac{dK^1}{dK^2}(s,\omega)\leq 1,\quad dK^2_s
(\omega)P(d\omega)-a.e. \,\,\mbox{on}\,\, [0,T]\times \Omega.
$$
\end{enumerate}
\end{definition}

We now introduce the following data :\\

$\bullet$ $L:=\left\{ L_{t},\,0\leq t\leq T\right\}$ and $U:=\left\{
U_{t},\,0\leq t\leq T\right\}$ are two real valued barriers
 which are $\cal P$-measurable and
continuous processes such that $L_t \leq U_t,\,\,\forall t\in[0,T]$.

$\bullet$ $\xi$ is an ${\cal F}_T$-measurable one dimensional random
variable such that $L_T\leq \xi\leq U_T$.

 $\bullet$ $f :
\Omega \times \left[ 0,T\right] \times \R^{1+d}\longrightarrow \R$
is a function which to $(t,\omega,y,z)$ associates $f(t,\omega,y,z)$
which is continuous with respect to $(y,z)$ and $\cal P$-measurable.

 $\bullet$ $g :
\Omega \times \left[ 0,T\right] \times \R\longrightarrow \R$ is a
function which to $(t,\omega,y)$ associates $g(t,\omega,y)$ which is
continuous with respect to $y$ and $\cal P$-measurable.\\

Let us now introduce the definition of our GBSDE with two reflecting
obstacles $L$ and $U$.
\begin{definition} We call $(Y,Z,K^+,K^-):=( Y_{t},Z_{t},K_{t}^+,K_{t}^-)_{t\leq T}$
a solution of the GBSDE with two reflecting barriers, associated
with coefficient $f ds +gdA_s+ dR_s$; terminal value $\xi$ and
barriers $L$ and $U$, if the following hold :
\begin{equation}
\label{eq0} \left\{
\begin{array}{ll}
(i) & 
 Y_{t}=\xi
+\integ{t}{T}f(s,Y_{s},Z_{s})ds+\dint_t^T dR_s+\dint_t^Tg(s,
Y_s)dA_s\\
&\qquad\quad+\integ{t}{T}dK_{s}^+ -\integ{t}{T}dK_{s}^-
-\integ{t}{T}Z_{s}dB_{s}\,, t\leq T,
\\ (ii)& Y\mbox{ between } L \mbox{ and } U,\, \, i.e.\,\,
\forall  t\leq T,\,\, L_t \leq Y_{t}\leq U_{t},\\ (iii) &\mbox{ the
Skorohod conditions hold : }\\ & \integ{0}{T}( Y_{t}-L_{t})
dK_{t}^+= \integ{0}{T}( U_{t}-Y_{t}) dK_{t}^-=0,\,\, \mbox{a.s.}, \\
(iv)& Y\in {\cal C} \quad K^+, K^-\in {\cal K} \quad Z\in {\cal
L}^{2,d},  \\ (v)& dK^+\perp  dK^-.
\end{array}
\right. \end{equation}
\end{definition}
\begin{remark}It should be pointed out here that in our setting we do not
require a square integrability on the solutions since no
integrability conditions are assumed on the data. This is not a
handicap since in many applications, such as stochastic games or
mathematical finance, we do not need such properties for the
solutions.
\end{remark}

 Next, we are going to suppose weaker conditions on the data under
which the GBSDE (\ref{eq0}) has a solution. We shall need the
following assumptions on $f$ and $g$ :
\medskip

\ni $(\bf{A.1})$ There exist two processes $\eta \in L^0(\Omega,
L^1([0,T], \R_+))$ and $C\in {\cal C}$ such that:
 $$\forall (s,\omega),\,\, |f(s, \omega, y, z )| \leq
 \eta_s(\omega)+\frac{C_s(\omega)}{2}|z|^2,\,\,
\forall y\in [L_s(\omega), U_s(\omega)],\,\, \forall z\in \R^d.$$

\ni $(\bf{A.2})$
 For all $(s,\omega),\,\, |g(s, \omega, y)| \leq 1,\,\, \forall y\in [L_s(\omega), U_s(\omega)].$

\medskip

For instance, Equation (\ref{eq0}) may not have a solution. Take,
for example, $L= U$ with $L$ not being a semi-martingale then
obviously we can not find a four-tuple which satisfies $ii)$ of
Equation (\ref{eq0}). Therefore, in order to obtain a solution, we
are led to assume :\\

\noindent$(\bf{A.3})$ There exists a continuous semimartingale $S_.
= S_0 + V_.^+ -V_.^- +\dint_0^. \alpha_s dB_s$, with $S_0\in\R, V^+,
V^-\in \cal K$ and $\alpha\in{\cal L}^{2,d}$, such that $$L_t \leq
S_t\leq U_t, \,\, \forall t\in [0,T].$$ \noindent$(\bf{A.4})$
$L_t\leq 0\leq U_t, \,\, \forall t\in [0,T].$
\\

Let us now give some remarks on the assumptions which show that our
assumptions are not restrictive and weaker than that on the current
literature.

\begin{remark}\label{rem1}
\begin{enumerate}
\item It is not difficult to see that if $L$ or $U$ is a continuous
semimartingale, then $(\bf{A.3})$ holds. Moreover if the barriers
processes $L$ and $U$ are completely separated on $[0,T]$, \it{i.e.}
$L_t < U_t,\,\, \forall t\in [0,T]$, then $(\bf{A.3})$ holds also
true. Indeed, let $\beta_t = \displaystyle\sup_{s\leq t} (\mid
L_s\mid +\mid U_s\mid )$. Since $L$ and $U$ are continuous then
$\forall t\in [0,T],$ $\bigg|\dfrac{U_t}{\beta_t}\bigg|\leq 1$ and
$\bigg|\dfrac{L_t}{\beta_t}\bigg|\leq 1$. It follows then from the
work \cite{HH} that there exists a continuous semimartingale
$\overline{S}$ such that $\dfrac{L_t}{\beta_t}\leq
\overline{S}_t\leq \dfrac{U_t}{\beta_t},\,\, \forall t\in [0,T]$.
Hence, the continuous semimartingale $\overline{S}\beta$ is between
$L$ and $U$.
\item By taking $Y_. -S_.$  instead of $Y_.$ one can
suppose, without loss of generality, that the semimartingale $S =
0$. Hence assumption $(\bf{A.4})$ will be assumed instead of
$(\bf{A.3})$.
\item Suppose that there exist $\widetilde{\eta}
\in L^0(\Omega, L^1([0,T], ds, \R_+))$ and $\widehat{\eta} \in
L^0(\Omega, L^1([0,T], dA_s, \R_+))$ such that : $\forall
(s,\omega)\in [0, T]\times\Omega,\,\,\forall y\in [L_s(\omega),
U_s(\omega)],\,\, \forall z\in \R^d,$
$$|f(s, \omega, y, z )| \leq
\widetilde{\eta}_s(\omega)+\phi(s, \omega, y) + \psi(s, \omega,
y)|z|^2,$$ and
 \begin{equation}\label{equa1}
 |g(s, \omega, y)| \leq
\widehat{\eta}_s(\omega)+\varphi(s, \omega, y),
\end{equation}
where $\phi$, $\psi$ and $\varphi$ are continuous functions on
$[0,T]\times \R$ and progressively measurable. Then conditions
$\bf{(A.1)}$ and $\bf{(A.2)}$ hold. In fact we just take, in
condition $(\bf{A.1})$, $\eta$ and $C$ as follows : $$
\begin{array}{lll}
& \eta_t(\omega) = \widetilde{\eta}_t(\omega)
+\displaystyle\sup_{s\leq t}\displaystyle\sup_{\alpha\in
[0,1]}|\phi(s,\omega, \alpha L_s +(1-\alpha) U_s)|,
\\ &  C_t(\omega) = 2\displaystyle\sup_{s\leq
t}\displaystyle\sup_{\alpha\in [0,1]}|\psi(s,\omega, \alpha L_s
+(1-\alpha) U_s)|.
\end{array}
$$ This means that the function $f$ can have, in particular, a
general growth in $y$ and quadratic growth in $z$. Now suppose that
the driver $g$ satisfies condition $(\ref{equa1})$, then for all $
(t,\omega)$ we have
$$
|g(t, \omega, y)| \leq \widehat{\eta}_t+\displaystyle\sup_{s\leq
t}\displaystyle\sup_{\alpha\in [0,1]}|\varphi(s,\omega, \alpha L_s
+(1-\alpha_s) U_s)|: = \overline{\eta}_t(\omega)\leq
\overline{\eta}_t(\omega) +1.
$$
 Now, if you take $\dfrac{g(t,y)}{1+\overline{\eta}_t}$ and
$(1+\overline{\eta}_t) dA_t$ instead of $g(t,y)$ and $dA_t$
respectively in equation (\ref{eq0})$(i)$, then we have condition
$(\bf{A.2})$.
\end{enumerate}
\end{remark}

\subsection{Exponential change for GRBSDE} The main idea for proving the existence of
 a solution for GRBSDE (\ref{eq0}) with data $(fds+gdA_s+dR_s,\xi, L, U)$ is to find a solution for a GRBSDE
with data obtained by improving an exponential transform of the data
$(fds+gdA_s+dR_s,\xi, L, U)$ which can be traced back to \cite{Kob}.
This transformation allows us, in particular, to bound the terminal
condition and the barriers associated with the transformed GRBSDE.
For this purpose, let us denote $|R|$ the total variation of the
process $R$ and define the processes $m, \overline{\xi},
\overline{L}, \overline{U}, \overline{g}, \overline{f}$,
$\overline{A}$ and $\overline{R}$ as follows:
$$
\begin{array}{lll}
&\ni \bullet\quad m_s =\displaystyle\sup_{r\leq s}|U_r|
+2\displaystyle\sup_{r\leq
 s}|C_r| +|R|_s + A_s
+1.
 \\ &\ni \bullet\quad\overline{\xi}= e^{ m_T(\xi-m_T)},\,\,
\overline{L}_s = e^{ m_s(L_s-m_s)},\quad \overline{U}_s = e^{
m_s(U_s-m_s)},
\\
& \\ &\bullet\quad \overline{g}(s,\overline{y})
=\dfrac{\widetilde{g}(s,(\overline{y}\vee
\overline{L}_s)\wedge\overline{U}_s)-4m_s}{8m_s},\quad\mbox{with}\quad\\
&\quad\,\, \widetilde{g}(s,\overline{y}) = \overline{y} \bigg(m_s
g(s, \dfrac{\ln(\overline{y})}{m_s}+m_s)\dfrac{dA_s}{dm_s}
+m_s\dfrac{dR_s}{dm_s} + (m_s-\dfrac{\ln(\overline{y})}{m_s})\bigg),
\,\, \overline{y}>0,
\\ & \\ & \bullet \quad\overline{f}(s,\overline{y},\overline{z})
=\widetilde{f}(s,(\overline{y}\vee
\overline{L}_s)\wedge\overline{U}_s,\overline{z})
-{\eta}_sm_s,\quad\mbox{with}\quad \\
&\quad\,\,\widetilde{f}(s,\overline{y},\overline{z}) = \overline{y}
\bigg(m_sf(s, \dfrac{\ln(\overline{y})}{m_s}+m_s,
\dfrac{\overline{z}}{m_s \overline{y}})
-\dfrac{|\overline{z}|^2}{2\overline{y}^2}\bigg),\quad
\overline{y}>0,\, \overline{z}\in\R^d,
\\ & \\ &\bullet \quad d\overline{A}_s = 8m_s dm_s\,\,\,\mbox{and}\,\,\,d\overline{R}_s =
\frac12 d\overline{A}_s +{\eta}_sm_s ds.
\end{array}
$$
%
\begin{remark} We should note here that $m$ is ${\cal
F}_t$-adapted, continuous and nondecreasing and then define a finite
variation process. This property will be used below.
\end{remark}


Suppose now that Equations $(\ref{eq0})$ has a solution $(Y, Z, K^+,
K^-)$ and define the processes $\overline{Y}, \overline{Z},
\overline{K}^+$ and $\overline{K}^-$ as follows :
\begin{equation}\label{eqq1}
\overline{Y}_. = e^{m_.(Y_.-m_.)},\quad\overline{Z}_. = m_.
\overline{Y}_.Z_.,\quad d\overline{K}_.^+ = m_. \overline{L}_.
dK^+_., \quad d\overline{K}_.^-  = m_. \overline{U}_. dK^-_..
\end{equation}
 Then
$(\overline{Y}, \overline{Z}, \overline{K}^+, \overline{K}^-)$ is
satisfying the following GBSDE
\begin{equation}
\label{eq2} \left\{
\begin{array}{ll}
(i) & 
 \overline{Y}_{t}=\overline{\xi} +\integ{t}{T}\overline{f}(s,\overline{Y}_{s},\overline{Z}_{s})ds +
 \integ{t}{T}\overline{g}(s,\overline{Y}_{s})d\overline{A}_s + \integ{t}{T} d\overline{R}_s\\
&\qquad\quad+\integ{t}{T}d\overline{K}_{s}^+
-\integ{t}{T}d\overline{K}_s^-
-\integ{t}{T}\overline{Z}_{s}dB_{s}\,, t\leq T,
\\ (ii)& \,\,
\forall  t\leq T,\,\, \overline{L}_t \leq \overline{Y}_{t}\leq
\overline{U}_{t},\\ (iii) &\integ{0}{T}(
\overline{Y}_{t}-\overline{L}_{t}) d\overline{K}_{t}^+=
\integ{0}{T}( \overline{U}_{t}-\overline{Y}_{t})
d\overline{K}_{t}^-=0,\,\, \mbox{a.s.} \\ (iv)& \overline{Y}\in
{\cal C} \quad \overline{K}^+, \overline{K}^-\in {\cal K} \quad
\overline{Z}\in {\cal L}^{2,d},
 \\ (v)&
d\overline{K}^+\perp d\overline{K}^-,
\end{array}
\right. \end{equation} where $\overline{\xi}, \overline{f},
\overline{g}, \overline{R}$, $\overline{L}$ and $\overline{U}$ are
given above.\\

More precisely, we have the following.
\begin{proposition} \label{pro3}
 Equations $(\ref{eq0})$ and $(\ref{eq2})$ are equivalent, in
the sense that if there exists a solution (resp. maximal solution)
to one of them then there exists a solution (resp. maximal solution)
for the other.
\end{proposition}
\bop. Suppose that Equation (\ref{eq0}) has a solution (resp.
maximal solution), say $(Y, Z, K^+, K^-)$. It follows from It\^{o}'s
formula that
$$
\begin{array}{lll}
 & e^{ m_t(Y_t-m_t)} \\ & = e^{m_T(Y_T-m_T)}+ \dint_{t}^{T}m_s
e^{m_s(Y_s-m_s)}(f(s,Y_{s},Z_{s})ds+ dK_s^+ -dK_s^- -Z_s dB_s)\\
& + \dint_{t}^{T}m_s e^{m_s(Y_s-m_s)}(g(s, Y_s)dA_s
+dR_s)+2\dint_{t}^{T}e^{m_s(Y_s-m_s)}m_s dm_s\\
&-\dint_{t}^{T}e^{m_sY_s}(Y_s-m_s)dm_s
 -\dfrac{1}{2}\dint_{t}^{T}e^{m_s(Y_s-m_s)}\mid m_sZ_{s}\mid^2
ds.
\end{array}
$$ Henceforth $$
\begin{array}{lll}
& e^{ m_t(Y_t-m_t)} \\ & = e^{m_T(\xi-m_T)}\ + \dint_{t}^{T}
e^{m_s(Y_s-m_s)}(m_s f(s,Y_{s},Z_{s}) -\dfrac{1}{2}\mid
m_sZ_{s}\mid^2)ds\\ & -\dint_{t}^{T} e^{m_s(Y_s-m_s)}m_sZ_s dB_s
+\dint_{t}^{T}m_s e^{m_s(Y_s-m_s)} dK_s^+
\\& -\dint_{t}^{T}m_s e^{m_s(Y_s-m_s)} dK_s^-
+\dint_{t}^{T}e^{m_s(Y_s-m_s)}(2m_s -Y_s)dm_s
\\ & +
\dint_{t}^{T}\bigg(e^{m_s(Y_s-m_s)}m_s g(s,
Y_s)\dfrac{dA_s}{dm_s}\bigg)dm_s
+\dint_{t}^{T}\bigg(e^{m_s(Y_s-m_s)} m_s
\dfrac{dR_s}{dm_s}\bigg)dm_s.
\end{array}
$$ Then it is clear that
$(\overline{Y}, \overline{Z}, \overline{K}^+, \overline{K}^-)$,
defined by (\ref{eqq1}) and  associated with coefficient
$\overline{f }ds +\overline{g}d\overline{A}_s+ d\overline{R}_s$, is
a solution (resp. maximal solution) of Equation $(\ref{eq2})$.
Conversely, Suppose that there exists a solution (resp. maximal
solution)
$(\overline{Y},\overline{Z},\overline{K}^+,\overline{K}^-)$ for
Equation $(\ref{eq2})$. Hence, by setting, for all $t\leq T$
$$ Y_t = \dfrac{\ln(\overline{Y}_t)}{m_t} +m_t, \quad Z_t =
\dfrac{\overline{Z}_t}{m_t {\overline{Y}_t  }}, \quad dK_t^{\pm} =
\dfrac{d\overline{K}_t^{\pm}}{m_t {\overline{Y}_t}},$$ one can see
that $(Y,Z,K^+,K^-)$ is a solution (resp. maximal solution) for
Equation $(\ref{eq0})$. \eop

The following proposition states some properties on the data
$(\overline{f} ds +\overline{g}d\overline{A}_s+ d\overline{R}_s,
\overline{\xi}, \overline{L}, \overline{U})$ of the transformed
GRBSDE (\ref{eq2}).

\begin{proposition}\label{pro0} Assume that assumptions $(\bf{A.1}),$ $(\bf{A.2})$ and $(\bf{A.4})$
hold. Then we have :
\begin{enumerate}
\item $\forall t\in [0,T],\,\,\, 0< \overline{L}_t \leq
e^{-m_t^2}\leq\overline{U}_t \leq
 e^{-1}< 1$ and $\overline{L}_T\leq\overline{\xi}\leq \overline{U}_T$.
\item The function $\widetilde{f}$ is $\cal P$-measurable and
continuous with respect to $(y,z)$ satisfying: $\forall
(s,\omega)\in [0,T]\times \Omega, \forall \overline{y}\in
[\overline{L}_s(\omega), \overline{U}_s(\omega)], \forall
\overline{z}\in \R^d$,
\begin{equation}\label{E1} -{\eta}_sm_s
-\dfrac{|\overline{z}|^2}{\overline{L}_s}\leq
\widetilde{f}(s,\omega, \overline{y},\overline{z}) \leq {\eta}_sm_s.
\end{equation}
\item For all $s\in [0,T],$ $\overline{y}\in \R$ and $\overline{z}\in
\R^d $
\begin{equation}\label{E2} -2{\eta}_sm_s -\dfrac{|\overline{z}|^2}{\overline{L}_s}\leq \overline{f}(s,
\overline{y},\overline{z}) \leq 0.
\end{equation}
\item For all $(s,\omega)\in [0,T]\times\Omega$ and $\overline{y}\in
[\overline{L}_s(\omega), \overline{U}_s(\omega)]$
\begin{equation}\label{eq19}
|\widetilde{g}(s, \overline{y})|\leq 4m_s,\quad \mbox{and} \quad
-1\leq \overline{g}(s, \overline{y})\leq 0.
\end{equation}
 \item $d\overline{R}$ is a positive measure.
\end{enumerate}
 \end{proposition}
\bop. Assertion $\it{1}.$ follows easily from assumption
$(\bf{A.4})$
and the fact that $m_t -1\geq U_t$,\,\,\, $\forall t\in [0,T]$.\\
Let us prove assertion $\it{2}$. It is not difficult to see that
$\widetilde{f}$ is $\cal P$-measurable and continuous with respect
to $(y,z)$ since $f$ is. It remains to prove inequality (\ref{E1}).
Let $(s,\omega)\in [0,T]\times\Omega, \,\,\overline{y}\in
[\overline{L}_s(\omega), \overline{U}_s(\omega)]$ and
$\overline{z}\in \R^d$, by condition $(\bf{A.1})$, we have
$$\begin{array}{ll}\widetilde{f}(s,\omega, \overline{y},\overline{z})
& \leq \overline{y}\bigg( m_s(\eta_s+ \dfrac{C_s}{2}
\dfrac{|\overline{z}|^2}{m_s^2 \overline{y}^2})
-\dfrac{|\overline{z}|^2}{2\overline{y}^2}\bigg)\\ &\leq
e^{-1}m_s\eta_s +
(\frac{C_s}{2m_s}-\frac{1}{2})\dfrac{|\overline{z}|^2}{\overline{y}}\\
& \leq m_s \eta_s,
\end{array}$$
since $\overline{y} \leq \overline{U}_s \leq e^{-1}<1$ and
$\frac{C_s}{2m_s}-\frac{1}{2} \leq 0$, $\quad \forall s\in [0,T]$ .
On the other hand, by using condition $(\bf{A.1})$, we get also that
$$\begin{array}{ll}\widetilde{f}(s,\omega, \overline{y},\overline{z})
& \geq \overline{y}\bigg( m_s(-\eta_s- \dfrac{C_s}{2}
\dfrac{|\overline{z}|^2}{m_s^2 \overline{y}^2})
-\dfrac{|\overline{z}|^2}{2\overline{y}^2}\bigg)\\ &\geq
-e^{-1}m_s\eta_s -
(\frac{C_s}{2m_s}+\frac{1}{2})\dfrac{|\overline{z}|^2}{\overline{y}}\\
& \geq - m_s\eta_s- \dfrac{|\overline{z}|^2}{\overline{L}_s},
\end{array}$$ since $\overline{y} \geq \overline{L}_s >0$ and $\frac{C_s}{2m_s}+\frac{1}{2} \leq 1$,$\quad \forall s\in
[0,T]$. Inequality (\ref{E1}) is then proved.
\\\\
Inequality (\ref{E2}) follows easily from inequality (\ref{E1}) and
then assertion \it{3}. holds.
\\ Assertions $\it{4}.$ and $\it{5}.$ follow immediately from
assumption $(\bf{A.2})$ and the definition of $m$.\eop

\begin{remark}
We should note here that, by taking advantage of Propositions 2.1
and 2.2, our problem is then reduced to find a maximal solution to
the following GRBSDE :
\begin{equation}
\label{eq221} \left\{
\begin{array}{ll}
(i) & 
 {Y}_{t}={\xi} +\integ{t}{T}{f}(s,{Y}_{s},{Z}_{s})ds +
 \integ{t}{T}g(s,{Y}_{s})d{A}_s + \integ{t}{T} d{R}_s\\
&\qquad\quad+\integ{t}{T}d{K}_{s}^+ -\integ{t}{T}d{K}_s^-
-\integ{t}{T}{Z}_{s}dB_{s}\,, t\leq T,
\\ (ii)&
\forall  t\leq T,\,\, {L}_t \leq {Y}_{t}\leq {U}_{t},\\ (iii)
&\integ{0}{T}({Y}_{t}-{L}_{t}) d{K}_{t}^+= \integ{0}{T}(
{U}_{t}-{Y}_{t}) d{K}_{t}^-=0,\,\, \mbox{a.s.}
\\ (iv)& {Y}\in {\cal C} \quad {K}^+,
{K}^-\in {\cal K} \quad {Z}\in {\cal L}^{2,d},
 \\ (v)&
d{K}^+\perp d{K}^-,
\end{array}
\right.
\end{equation} under the following assumptions :\\\\
\ni $(\bf{H.0})$ $dR \geq 0$, \it{i.e.} $R\in \cal K$.

\ni $(\bf{H.1})$ There exist two processes $\eta \in L^0(\Omega,
L^1([0,T], \R_+))$ and $C\in {\cal C}$ such that:
 $$\forall (s,\omega),\quad -\eta_s(\omega)-\frac{C_s(\omega)}{2}|z|^2\leq f(s, \omega, y, z ) \leq 0,\quad
\forall y\in \R,\quad \forall z\in \R^d.$$

\ni $(\bf{H.2})$
 $\forall (s,\omega),\quad -1\leq g(s, \omega, y) \leq 0,\quad
\forall y\in \R$.
\\\\
\noindent$(\bf{H.3})$ $0<L_t\leq U_t< 1,  \,\, \forall t\in [0,T].$

\noindent$(\bf{H.4})$ There exists a continuous  nondecreasing
process $S = S_0 -V $, with $S_0\in\R, \,\, V\in \cal K$, such that
$L_t \leq S_t\leq U_t, \,\, \forall t\in [0,T].$
\end{remark}

We devote the next section to the existence of maximal solution for
GBSDE (\ref{eq221}) under assumptions $(\bf{H.0})-(\bf{H.4})$ and
then equivalently  to the existence of maximal solution for GBSDE
(\ref{eq0}).
\section{Existence of maximal solution for GBSDE (\ref{eq221})}
Our objective now is to prove that, under assumptions
$(\bf{H.0})$--$(\bf{H.4})$, GRBSDE (\ref{eq221}) has a maximal
solution $(Y_t ,Z_t ,K^+_t , K_t^- )_{t\leq T}$, in the sense that
for any other $(Y_t^{'} ,Z_t^{'} ,K^{'+}_t , K_t^{'-} )_{t\leq T}$
of (\ref{eq221}) we have for all $t \leq T$, $Y_t \geq Y_t^{'}$,
$P$-a.s. The proof of our result is based on regularization by
sup-convolution techniques and a truncation procedure by means of a
family of stopping times.
\subsection{Approximations}

It is not difficult to prove the following lemma which gives an
approximation of continuous functions by Lipschitz functions (see,
for example, Lepeltier and San Martin \cite{LS}).
\begin{lemma}\label{lem0} Let $f_n$ and $g_n$ be two sequences
of functions defined by
\begin{equation}\label{eqq2}
 f_n(t,y,z) = \sup_{p\in \R, q\in \R^d} \{ {f}(t, p, q)\vee
(-n)-n|p-y| - n|q-z| \},
\end{equation}
 and
\begin{equation}\label{eqq3}
g_n(t,y) = \sup_{p\in \R} \{{g}(t, p)\vee (-n)-n|p-y| \}.
\end{equation} Assume that assumptions $(\bf{H.0})$--$(\bf{H.4})$ hold. Then we have the following :
\begin{enumerate}
\item For all $(t,\omega, y,z,n)\in [0,T]\times\Omega\times\R\times\R^d\times\N,\,\,$
$$f_0(t,y,z)=0\geq f_n(t,y,z)
\geq f_{n+1}(t,y,z)\geq {f}(t, y, z) \geq -{\eta}_t
-\dfrac{C_t}{2}|z|^2.$$
\item For all $(t,\omega, y,n)\in [0,T]\times\Omega\times\R\times\N,\,\,$
$$g_0(t,y)=0\geq g_n(t,y) \geq g_{n+1}(t,y)\geq {g}(t, y) \geq -1.$$
\item For all $(t,\omega,y, z, n)\in [0,T]\times\Omega\times\R\times\R^d\times\N,$
$$
-n\leq f_n(s,y,z)\leq 0.
$$
\item $f_n$ is uniformly $n$-Lipschitz with respect to $(y,z)$.
\item $g_n$ is uniformly $n$-Lipschitz with respect to $y$.
\item For all $(t,\omega)\in [0,T]\times\Omega$, $ (f_n(t,y,z)) _{n\geq 0}$ converges to ${f}(t, y, z)$
as $n$ goes to $+\infty$ uniformly on every compact of $\R\times
\R^d.$
\item For all $(t,\omega)\in [0,T]\times\Omega$, $ (g_n(t,y)) _{n\geq 0}$ converges
to ${g}(t, y)$ as $n$ goes to $+\infty$ uniformly on every compact
of $\R$.
\end{enumerate}
\end{lemma}
 Since the integrability conditions on the data of our GRBSDE are weaker, it will be useful to define the following family of
stopping times $(\tau_n)_{n\geq 0}$
\begin{equation}\label{tau}
\tau_n =\inf\{s\geq 0 : A_s+R_s +C_s+\dint_0^s\eta_r dr\geq
n\}\wedge T.
\end{equation}
\begin{remark} We should note here that the family $(\tau_n)_{n\geq 0}$
satisfies the following property which will be useful in the sequel
$$P\bigg[\displaystyle{\bigcup_{j=1}^{\infty}(\tau_j= T)}\bigg]
=1.$$ We say that $(\tau_n)_{n\geq 0}$ is a stationary sequence of
stopping times.

\ni Indeed, let
$\displaystyle{\omega\in\bigcap_{j=1}^{\infty}(\tau_j< T)}$ then
$\forall j\geq 1, \,\, M_T(\omega):= (A_T+R_T +C_T+\dint_0^T\eta_r
dr)(\omega)\geq j$ and hence $M_T(\omega) =+\infty$. Therefore
$P\bigg[\displaystyle\bigcap_{j=1}^{\infty}(\tau_j< T)\bigg]\leq
P[M_T = +\infty] =0$ and then
$P\bigg[\displaystyle{\bigcup_{j=1}^{\infty}(\tau_j= T)}\bigg]
=1$.\eop
\end{remark}
 Set $dA_s^n = 1_{\{s\leq \tau_n\}}dA_s, n\in \N$  and
$dR_s^i = 1_{\{s\leq \tau_i\}}dR_s, i\in\N$ and consider the
following BSDE with two reflecting barriers
\begin{equation}
\label{eq3} \left\{
\begin{array}{ll}
(i) & 
 Y_{t}^{n, i}=\xi +
\integ{t}{T}f_n(s,Y_{s}^{n, i},Z_{s}^{n, i})ds+ \integ{t}{T} g_n(s,Y_s^{n, i})dA_s^n + \integ{t}{T}dR_s^i \\
&\qquad\quad+\integ{t}{T}dK_{s}^{n,i+} -\integ{t}{T}dK_{s}^{n,i-}
-\integ{t}{T}Z_{s}^{n, i}dB_{s}\,, t\leq T,
\\ (ii)& \,\,
\forall  t\leq T,\,\, {L}_t \leq Y_{t}^{n, i}\leq {U}_{t},\\
(iii) & \integ{0}{T}( Y_{t}^{n, i}-{L}_{t}) dK_{t}^{n,i+}=
\integ{0}{T}(
{U}_{t}-Y_{t}^{n, i}) dK_{t}^{n,i-}=0,\,\, \mbox{a.s.}\\
(iv)&Y^{n, i}\in {\cal C} \quad K^{n,i+}, K^{n,i-}\in {\cal K} \quad
Z^{n, i}\in {\cal L}^{2,d},  \\ (v)& d{K}^{n,i+}\perp d{K}^{n,i-}.
\end{array}
\right. \end{equation} It follows from Theorem \ref{th22} (see
Appendix) that Equation (\ref{eq3}) has a unique solution. Moreover,
for all $n$ and $i$
\begin{equation}\label{eqq4}
 \E\sup_{t\leq T}|Y_t^{n,i}|^2 + \E\dint_0^T \mid Z^{n, i}_s\mid^2 ds+\E(K_T^{n,i+})^2 <+\infty.
\end{equation}

 Our objective now is to study the GRBSDE (\ref{eq3}) and justify the
passage to the limit and identify the limit as the solution of our
equation.

\subsection{The study of Equation (\ref{eq3}) for $n$ fixed}
The following result follows easily from the Comparison theorem
(Theorem \ref{th13} in Appendix).
\begin{proposition}\label{pro1} Let us suppose that assumptions $(\bf{H.0})$--$(\bf{H.4})$
hold. Then we have the following.
\begin{itemize}
\item[i)] Fix $n$, we get for all $i\geq 0$ and $t\leq T$
$$
{L}_t\leq Y_{t}^{n, i}\leq Y_{t}^{n, i+1}\leq {U}_{t},\quad
dK^{n,i+}\geq dK^{n,i+1+} \quad\mbox{ and }\quad dK^{n, i+1-}\geq
dK^{n,i-}.
$$
\item[ii)] Fix $i$, we get for all $n\geq
0$ and $t\leq T$
$$
{L}_t\leq Y_{t}^{n+1, i}\leq Y_{t}^{n, i}\leq {U}_{t},\quad
dK^{n,i+}\leq dK^{n+1,i+} \quad\mbox{ and }\quad dK^{n+1, i-}\leq
dK^{n,i-}.
$$
\end{itemize}
\end{proposition}
\bop. Since the family of stopping times $(\tau_i)_{i\geq 0}$ is
increasing then, for all $i\geq 0$,  $dR^i\leq dR^{i+1}$. The
assertion \it{i)} follows by using inequality (\ref{eqq4}) and
Theorem \ref{th13} in Appendix.\\
Assertion \it{ii)} follows also by using Lemma \ref{lem0} and
Theorem \ref{th13} in Appendix.\eop

\ni Now we want to study Equation (\ref{eq3}) when
$n$ is fixed. Let us set\\
$\bullet$ $Y^n = \displaystyle\sup_{i}Y^{n, i}$
\\
$\bullet$ $dK^{n-} = \displaystyle\sup_{i}dK^{n,i-}$ which is a
positive measure. \\
$\bullet$ $dK^{n+} = \displaystyle\inf_{i}dK^{n, i+}$ which is also
a positive measure since $K_T^{n, 0+}<+\infty, \, P-a.s.$

\begin{proposition}\label{pro2}Assume that assumptions $(\bf{H.0})$--$(\bf{H.4})$
hold. Then we have the following.
\begin{enumerate}
\item There exists a process $Z^n\in {\cal L}^{2,d}$ such that, for
all $j\in\N$,
$$
\E\integ{0}{\tau_j}|Z_s^{n, i} - Z_s^{n}|^2ds \longrightarrow
0,\quad \mbox{as $i$ goes to infinity}.
$$
\item The process $(Y^n, Z^n, K^{n+}, K^{n-})$ is the unique solution
of the following GBSDE with two reflecting barriers
\begin{equation}
\label{eqq} \left\{
\begin{array}{ll}
(i) & 
 Y_{t}^{n}={\xi} +
\integ{t}{T}f_n(s,Y_{s}^{n},Z_{s}^{n})ds + \integ{t}{T} g_n(s,Y_s^{n})dA_s^n + \integ{t}{T}dR_s\\
&\qquad\quad+\integ{t}{T}dK_{s}^{n+} -\integ{t}{T}dK_{s}^{n-}
-\integ{t}{T}Z_{s}^{n}dB_{s}\,, t\leq T,
\\ (ii)&
\forall  t\leq T,\,\, {L}_t \leq Y_{t}^{n}\leq {U}_{t},\\
(iii) & \integ{0}{T}( Y_{t}^{n}-{L}_{t}) dK_{t}^{n+}= \integ{0}{T}(
{U}_{t}-Y_{t}^{n}) dK_{t}^{n-}=0,\,\, \mbox{a.s.}\\
(iv)&Y^{n}\in {\cal C} \quad K^{n+}, K^{n-}\in {\cal K} \quad
Z^{n}\in {\cal L}^{2,d},
\\ (v)&
dK^{n+}\perp dK^{n-}.
\end{array}
\right. \end{equation}
\end{enumerate}
\end{proposition}\bop.
\ni \it{1.} Let $j, i, i'\in\N$ such that $j\leq i\leq i'$ and $t\in
[0, \tau_j]$ where $\tau_j$ is defined in (\ref{tau}). Since the
family of stopping times $(\tau_j)_j$ is increasing, we have
$$
\integ{t}{\tau_j}(dR_s^i-dR_s^{i'})=\integ{t}{\tau_j\wedge\tau_i}dR_s-\integ{t}{\tau_j\wedge\tau_{i'}}dR_s=
\integ{t}{\tau_j}dR_s-\integ{t}{\tau_j}dR_s=0
$$
Moreover, by taking advantage of the fact that $f_n$ and $g_n$ are
$n-$Lipschitz, we get
$$
\begin{array}{ll}
& Y_{t}^{n, i} -Y_{t}^{n,i'} \\ &= Y_{\tau_j}^{n, i}
-Y_{\tau_j}^{n,i'} +\integ{t}{\tau_j}(f_n(s,Y_{s}^{n, i},Z_{s}^{n,
i}) -f_n(s,Y_{s}^{n,i'},Z_{s}^{n,i'}))ds\\ & + \integ{t}{\tau_j}
(g_n(s,Y_s^{n, i})-g_n(s,Y_s^{n, i'}))dA_s^n
+\integ{t}{\tau_j}(dK_{s}^{n,i+}-dK_{s}^{n,i'+})\\ &
-\integ{t}{\tau_j}(dK_{s}^{n,i-}-dK_{s}^{n,i'-})
-\integ{t}{\tau_j}(Z_{s}^{n, i}- Z_{s}^{n,i'})dB_{s}\\ &=
 Y_{\tau_j}^{n, i} -Y_{\tau_j}^{n,i'} +\integ{t}{\tau_j}\bigg(\alpha_s^{n,i,i'
}(Y_s^{n,i} -Y_{s}^{n,i'}) +\langle \beta_s^{n,i,i'
}, Z_s^{n,i} -Z_{s}^{n,i'}\rangle \bigg)ds \\
&+ \integ{t}{\tau_j} C_s^{n,i,i' }(Y_s^{n,i} -Y_{s}^{n,i'})
dA_s^n+\integ{t}{\tau_j}(dK_{s}^{n,i+}-dK_{s}^{n,i'+})\\ &
-\integ{t}{\tau_j}(dK_{s}^{n,i-}-dK_{s}^{n,i'-})
-\integ{t}{\tau_j}(Z_{s}^{n,i}- Z_{s}^{n,i'})dB_{s},
\end{array} $$
where $C^{n,i,i'}, \alpha^{n,i,i'}$ and $\beta^{n,i,i'}$ are bounded
by $n$ and $\langle . , .\rangle$ denotes the scalar product in $\R^d$.\\
Set $e_t^n:= e^{2nA_t^n +(2n^2 +2n )t}$. Applying It\^{o}'s formula
to $(Y_{t}^{n, i} -Y_{t}^{n,i'})^2 e_t^n$ we get
$$
\begin{array}{ll}
&(Y_{t}^{n, i} -Y_{t}^{n,i'})^2 e_t^n +
\integ{t}{\tau_j}e_s^n|Z_s^{n, i} - Z_s^{n,i'}|^2ds
\\ & \qquad+\integ{t}{\tau_j}e_s^n \bigg(2n(Y_{s}^{n, i} -Y_{s}^{n,i'})^2
dA_s^n +(2n^2 +2n )(Y_{s}^{n, i} -Y_{s}^{n,i'})^2 ds\bigg)\\ & =
(Y_{\tau_j}^{n, i} -Y_{\tau_j}^{n,i'})^2 e_{\tau_j}^n+
2\integ{t}{\tau_j}e_s^n(Y_{s}^{n, i} -Y_{s}^{n,i'})
\bigg(\alpha_s^{n,i,i' }(Y_s^{n,i} -Y_{s}^{n,i'}) +\langle
\beta_s^{n,i,i' }, Z_s^{n,i} -Z_{s}^{n,i'}\rangle\bigg)ds\\ &+
2\integ{t}{\tau_j}C_s^{n,i,i'}e_s^n(Y_{s}^{n, i} -Y_{s}^{n,i'})^2
dA_s^n + \integ{t}{\tau_j}e_s^n(Y_{s}^{n, i}
-Y_{s}^{n,i'})(dK_{s}^{n,i+}-dK_{s}^{n,i'+})
\\ &-\integ{t}{\tau_j}e_s^n(Y_{s}^{n, i} -Y_{s}^{n,i'})(dK_{s}^{n,i-}-dK_{s}^{n,i'-})
-2\integ{t}{\tau_j}e_s^n(Y_{s}^{n, i} -Y_{s}^{n,i'})( Z_{s}^{n,i}-
Z_{s}^{n,i'})dB_{s}.
\end{array} $$
In force of Proposition \ref{pro1}, we may conclude that for $t\in
[0, \tau_j]$,
$$
\integ{t}{\tau_j}e_s^n(Y_{s}^{n, i}
-Y_{s}^{n,i'})(dK_{s}^{n,i+}-dK_{s}^{n,i'+})\leq 0\,\,\mbox{and}\,\,
-\integ{t}{\tau_j}e_s^n(Y_{s}^{n, i}
-Y_{s}^{n,i'})(dK_{s}^{n,i-}-dK_{s}^{n,i'-})\leq 0.
$$
and
$$
\begin{array}{ll}
& 2\integ{t}{\tau_j}e_s^n(Y_{s}^{n, i} -Y_{s}^{n,i'})\langle
\beta_s^{n,i,i' }, Z_s^{n,i} -Z_{s}^{n,i'}\rangle ds
\\ &\leq  2n^2\integ{t}{\tau_j}e_s^n(Y_{s}^{n, i}
-Y_{s}^{n,i'})^2ds + \frac12 \integ{t}{\tau_j}e_s^n|Z_{s}^{n, i}
-Z_{s}^{n,i'}|^2ds.
\end{array}
$$
Hence it follows that
$$
\begin{array}{ll}
& (Y_{t}^{n, i} -Y_{t}^{n,i'})^2 e_t^n +\dfrac12
\integ{t}{\tau_j}e_s^n|Z_s^{n, i} - Z_s^{n,i'}|^2ds
\\ &\leq (Y_{\tau_j}^{n, i}
-Y_{\tau_j}^{n,i'})^2
e_{\tau_j}^{n}-2\integ{t}{\tau_j}e_s^n(Y_{s}^{n, i} -Y_{s}^{n,i'})
(Z_{s}^{n,i}- Z_{s}^{n,i'}) dB_{s} ,
\end{array}
$$
Since $\integ{0}{.\,\wedge \tau_j}e_s^n(Y_{s}^{n, i} -Y_{s}^{n,i'})
(Z_{s}^{n,i}- Z_{s}^{n,i'})dB_{s}$ is an $({\cal F}_{t}, P)$-
martingale we obtain
$$
\begin{array}{ll}
 \E\integ{0}{\tau_j}e_s^n|Z_s^{n, i} - Z_s^{n,i'}|^2ds
\leq 2\E e_{\tau_j}^{n}(Y_{\tau_j}^{n, i} -Y_{\tau_j}^{n,i'})^2,
\end{array}
$$
 According to Bulkholder-Davis-Gundy
inequality, there exists a constant $C >0$ such that
\begin{equation}\label{eqq6}
\E\bigg[\sup_{t\leq \tau_j }(Y_{s}^{n, i} -Y_{s}^{n,i'})^2 +\dfrac12
\integ{0}{\tau_j }|Z_s^{n, i} - Z_s^{n,i'}|^2ds\bigg] \leq 4 C
 \E e_{\tau_j}^{n}(Y_{\tau_j }^{n, i} -Y_{\tau_j }^{n,i'})^2.
\end{equation}
 By using that $Y^{n,i}$ are bounded by $1$, assumption $(\bf{H.3})$ and Lebesgue's
dominated convergence theorem we have $\forall j\in \N$,
$$
\E\bigg[Y_{\tau_j }^{n, i} -Y_{\tau_j
}^{n,i'}\bigg]^2\longrightarrow 0,\,\,\,\mbox{as}
\,\,\,i\,\,\,\mbox{ goes to}\,\,\, +\infty.
$$
As a consequence, there exists $Z^n \in {\cal L}^{2,d}$ such that
$\forall j\in \N$
$$ \E \integ{0}{\tau_j }|Z_s^{n, i}- Z_s^{n}|^2ds \longrightarrow
0,\,\,\,\mbox{as} \,\,\,i\,\,\,\mbox{ goes to}\,\,\, +\infty.$$
Assertion \it{1.} is then proved.\\\\
 \ni \it{2.} Let us first prove assertion $(iv)$ of Equation
(\ref{eqq}). In view of passing to the limit in Inequality
(\ref{eqq6}) we get
$$
\E\bigg[\sup_{t\leq \tau_j }(Y_{s}^{n, i} -Y_{s}^{n})^2\bigg]
\longrightarrow 0,\,\,\,\mbox{as} \,\,\,i\,\,\,\mbox{ goes to}\,\,\,
+\infty,
$$
and then we can conclude that $Y^n$ is continuous, \it{i.e.} $Y^{n}\in {\cal C}$.\\
  It is
clear, from Proposition \ref{pro1}, that $K^{n,i+}$ converges to the
continuous and increasing process $K^{n+}$. Moreover,
$\E(K_T^{n+})^2\leq \E(K_T^{n,0+})^2 <+\infty$, for all $n$.
Therefore $K^{n+}\in {\cal K}$.\\
Now, passing to the limit in Equation (\ref{eq3})$(i)$ on
$[0,\tau_j]$, we get also that $\E(K_{\tau_j}^{n-})^2 < +\infty,\,\,
\forall j, n$. Since $P[\displaystyle{\cup_{j=1}^{\infty}(\tau_j=
T)}] =1$, we get $K_{T}^{n-} < +\infty,\,\, \forall n\,\,\,P$-a.s.
Then $K^{n-}\in {\cal K}$. Consequently, assertion $(iv)$ of
Equation (\ref{eqq}) is proved.
\\
 Let us now show that $(Y^n, Z^n, K^{n+}, K^{n-})$ satisfies $(i)$. In view of passing to the limit, as $i$ goes to infinity, in the following
equation
$$
\begin{array}{ll}
 Y_{t}^{n, i}=& Y_{\tau_j}^{n, i} + \integ{t}{\tau_j}f_n(s,Y_{s}^{n,
i},Z_{s}^{n, i})ds+ \integ{t}{\tau_j} g_n(s,Y_s^{n, i})dA_s^n +
\integ{t}{\tau_j}dR_s^i \\ &+\integ{t}{\tau_j}dK_{s}^{n,i+}
-\integ{t}{\tau_j}dK_{s}^{n,i-} -\integ{t}{\tau_j}Z_{s}^{n,
i}dB_{s}.
\end{array}
$$
We obtain, $P-$a.s.
$$
\begin{array}{ll}
Y_{t}^{n}= & Y_{\tau_j}^{n} +
\integ{t}{\tau_j}f_n(s,Y_{s}^{n},Z_{s}^{n})ds+ \integ{t}{\tau_j}
g_n(s,Y_s^{n})dA_s^n + \integ{t}{\tau_j}dR_s
\\ & +\integ{t}{\tau_j}dK_{s}^{n,+} -\integ{t}{\tau_j}dK_{s}^{n,-}
-\integ{t}{\tau_j}Z_{s}^{n}dB_{s},
\end{array}
$$
Since $(\tau_j)_{j\geq 0}$ is a stationary sequence of stopping
times we get $P-$a.s.
$$
\begin{array}{ll}
Y_{t}^{n}= & \xi + \integ{t}{T}f_n(s,Y_{s}^{n},Z_{s}^{n})ds+
\integ{t}{T} g_n(s,Y_s^{n})dA_s^n + \integ{t}{T}dR_s
\\ & +\integ{t}{T}dK_{s}^{n,+} -\integ{t}{T}dK_{s}^{n,-}
-\integ{t}{T}Z_{s}^{n}dB_{s},
\end{array}
$$
Hence the process $(Y^n, Z^n, K^{n+}, K^{n-})$ satisfies $(i)$ of
Equation (\ref{eqq}). \\ We now prove that the Skorohod conditions
$(iii)$ of Equation (\ref{eqq}) is satisfied.\\
Since $Y^n = \displaystyle\sup_{i}Y^{n, i}$, we have
$$
0\leq \integ{0}{T}( {U}_{t}-Y_{t}^{n})
dK_{t}^{n,i-}\leq\integ{0}{T}( {U}_{t}-Y_{t}^{n,i}) dK_{t}^{n,i-}=0.
$$
Therefore
$$
\integ{0}{T}( {U}_{t}-Y_{t}^{n})dK_{t}^{n,i-} = 0.
$$
 It follows then from Fatou's lemma that
$$
0\leq\integ{0}{T}( {U}_{t}-Y_{t}^{n}) dK_{t}^{n-}\leq
\liminf_{i}\integ{0}{T}( {U}_{t}-Y_{t}^{n})dK_{t}^{n,i-} = 0.
$$
Consequently
$$
\integ{0}{T}( {U}_{t}-Y_{t}^{n}) dK_{t}^{n-}  = 0.
$$
On the other hand
$$
0\leq \integ{0}{T}(Y_{t}^{n, i} -L_t) dK_{t}^{n+} \leq
\integ{0}{T}(Y_{t}^{n, i} -L_t) dK_{t}^{n,i+} = 0.
$$
Hence
$$
\integ{0}{T}(Y_{t}^{n, i} -L_t) dK_{t}^{n+} = 0.
$$
Applying Fatou's lemma we obtain
$$
0\leq\integ{0}{T}(Y_{t}^{n} -L_t) dK_{t}^{n+} \leq \liminf_{i}
\integ{0}{T}(Y_{t}^{n, i} -L_t) dK_{t}^{n+}= 0.
$$
Henceforth
$$
\integ{0}{T}(Y_{t}^{n} -L_t) dK_{t}^{n+} = 0.
$$
Let $A^{n,i}\in{\cal P}$ such that $dK^{n,i+}$ is supported by
$A^{n,i}$ and $dK^{n,i-}$ is supported by $(A^{n,i})^c$, then
$$
0\leq \E\dint_0^T 1_{\bigcup_{i}(A^{n,i})^c}dK^{n+}\leq
\sum_{i}\E\dint_0^T 1_{(A^{n,i})^c}dK^{n+}\leq \sum_{i}\E\dint_0^T
1_{(A^{n,i})^c}dK^{n, i+} =0.
$$
It follows that $dK^{n+}$ is supported by
$\displaystyle\bigcap_{i}A^{n,i}$. Now
$$
\E\dint_0^T 1_{(\bigcap_{i}A^{n,i})}dK^{n-}\leq
\sup_j\bigg(\E\dint_0^T 1_{(\bigcap_{i}A^{n,i})}dK^{n,j-}\bigg)\leq
\sup_j\bigg(\E\dint_0^T 1_{A^{n,j}}dK^{n,j-}\bigg) =0.
$$
Hence, the measures $dK^{n+} = \displaystyle\inf_{i}dK^{n, i+}$,\,\,
$dK^{n-} = \displaystyle\sup_{i}dK^{n,i-}$ are also singular. The
proof of Proposition 3.2 is finished.
 \eop
 \subsection{The study of the Equation (\ref{eqq})} Let $(Y^n, Z^n, K^{n+}, K^{n-})$
be the process given in Proposition \ref{pro2} the unique solution
of the GBSDE with two reflecting barriers (\ref{eqq}).\\
We should recall here that the result of the previous section gives
us the following :
\begin{equation}\label{eqq9}
\forall j \in\N,\,\, \forall n\in\N, \quad E(K_T^{n+})^2
+\E\dint_0^{\tau_j} \mid Z_s^n \mid^2ds <+\infty
\end{equation}
and then $\E(K_{\tau_j}^{n-})^2 < +\infty$, for all $j$
and $n$.\\

In view of passing to the limit in Proposition \ref{pro1}, we get
the following.
\begin{proposition}
For all $n\geq 0$, we obtain
$$ Y^{n+1}\leq Y^n, \qquad dK_{t}^{n+}
\leq dK_{t}^{n+1+}, \qquad dK_{t}^{n+1-} \leq dK_{t}^{n-}.
$$

\end{proposition}
 In order to study Equation (\ref{eqq}), let us set
 \begin{equation}\label{eqq5}
 \begin{array}{lll}
 & \bullet\,\,\, Y = \displaystyle\inf_{n}Y^{n}.\\ &
\bullet\,\,\, dK^{+} = \displaystyle\sup_{n}dK^{n+},\,\,\mbox{ which
is also a positive measure}.
\\ &
\bullet\,\,\, dK^{-} = \displaystyle\inf_{n}dK^{n-},\,\,\mbox{ which
is also a positive measure since}\,\, K_T^{n-}<+\infty, \, P-a.s.
\end{array}
\end{equation}
The following result states the convergence of the process $Z^n$ in
$L^2([0,\tau_j]\times \Omega)$.
\begin{proposition}\label{th1} Assume that assumptions $(\bf{H.0})$--$(\bf{H.4})$
hold. Then there exists a process $Z\in {\cal L}^{2,d}$ such that,
for all $j$,
$$
\E\integ{0}{\tau_j}|Z_s^{n} - Z_s|^2ds \longrightarrow 0,\quad
\mbox{as $n$ goes to infinity}.
$$
\end{proposition}
\bop. For $s\in [0, 1]$ and $j\in \N$, let us set $\psi(s) =
\dfrac{e^{12j s}-1}{12j} -s$. We mention that $\psi$ satisfies the
following for all  $s\in [0,1]$,
\begin{equation} \label{eqq7}
\begin{array}{lll}
 & \psi'(s) = e^{12j s} -1 ,\,\,\,  \psi''(s) = 12j e^{12j s} = 12j
 \psi'(s)+12j
 \\ &
0\leq 12j s\leq \psi'(s)\leq 12 j e^{12j }s\leq 12 j e^{12j }.
\end{array}
\end{equation}
Recall that $dA_s^n = 1_{\{s\leq \tau_n\}}dA_s,$ for $n\in \N$.
Applying It\^{o}'s formula to $\psi(Y^n -Y^m)$, we get for $m\geq n$
and $t\leq \tau_j$,
$$
\begin{array}{ll}
&\psi(Y_{t}^n -Y_{t}^{m}) \\&= \psi(Y_{\tau_j}^n -Y_{\tau_j}^{m})
-\integ{t}{\tau_j}
g_m(s,Y_{s}^{m})\psi'(Y_{s}^n -Y_{s}^{m})1_{\{\tau_n \leq s\leq \tau_m\}} dA_s \\
&+\integ{t}{\tau_j} \bigg(g_n(s,Y_s^n)
-g_m(s,Y_{s}^{m})\bigg)\psi'(Y_{s}^n -Y_{s}^{m})1_{\{s\leq \tau_n\}}
dA_s\\ & +\integ{t}{\tau_j}(f_n(s,Y_{s}^n,Z_{s}^n)
-f_m(s,Y_{s}^{m},Z_{s}^{m}))\psi'(Y_{s}^n -Y_{s}^{m})ds
\\ &
+\integ{t}{\tau_j}\psi'(Y_{s}^n -Y_{s}^{m}) d(K_{s}^{n+}-K_{s}^{m+})
-\integ{t}{\tau_j}\psi'(Y_{s}^{n} -Y_{s}^{m}) d(K_{s}^{n-}-K_{s}^{m-}) \\
&-\integ{t}{\tau_j}\psi'(Y_{s}^n -Y_{s}^{m})(Z_{s}^n-
Z_{s}^{m})dB_{s} -\dfrac12\integ{t}{\tau_j}\psi''(Y_{s}^n
-Y_{s}^{m})|Z_{s}^n- Z_{s}^{m}|^2 ds.
\end{array}
$$
Since, $K^{n+}$ (resp. $K^{m+}$) moves only when $Y^n$ (resp. $Y^m$)
reaches the obstacles ${L}$, $Y^m \leq Y^n$ and $\psi'(0) =0$, we
have $$ \integ{t}{\tau_j}\psi'(Y_{s}^n -Y^m_s)
d(K_{s}^{n+}-K_{s}^{m+}) = -\integ{t}{\tau_j}\psi'(Y_{s}^n -{L}_{s})
dK_{s}^{m+}.$$ By the same way we get also that
 $$
\integ{t}{\tau_j}\psi'(Y_{s}^n -Y^m_s) d(K_{s}^{n-}-K_{s}^{m-}) =
\integ{t}{\tau_j}\psi'({U}_{s}-Y_{s}^m) dK_{s}^{n-}.$$ Henceforth
$$
\begin{array}{ll}
& \psi(Y_{t}^n -Y_{t}^{m}) \\ &= \psi(Y_{\tau_j}^n -Y_{\tau_j}^{m})
- \integ{t}{\tau_j}g_m(s,Y_{s}^{m})\psi'(Y_{s}^n
-Y_{s}^{m})1_{\{\tau_n \leq s\leq \tau_m\}} dA_s\\
&+\integ{t}{\tau_j}
\bigg(g_n(s,Y_{s}^{n})-g_m(s,Y_{s}^{m})\bigg)\psi'(Y_{s}^n
-Y_{s}^{m})1_{\{s\leq \tau_n\}}dA_s\\ &
+\integ{t}{\tau_j}\bigg(f_n(s,Y_{s}^n,Z_{s}^n)
-f_m(s,Y_{s}^{m},Z_{s}^{m})\bigg)\psi'(Y_{s}^n -Y_{s}^{m})ds
\\ &
-\integ{t}{\tau_j}\psi'(Y_{s}^n -{L}_{s}) dK_{s}^{m+}-
\integ{t}{\tau_j}\psi'({U}_{s}-Y_{s}^m) dK_{s}^{n-}\\
&-\integ{t}{\tau_j}\psi'(Y_{s}^n -Y_{s}^{m})(Z_{s}^n-
Z_{s}^{m})dB_{s} -\dfrac12\integ{t}{\tau_j}\psi''(Y_{s}^n
-Y_{s}^{m})|Z_{s}^n- Z_{s}^{m}|^2 ds.
\end{array}
$$ In view of Lemma \ref{lem0} (assertions \it{1} and \it{2}) and Inequality (\ref{eqq7}), we get
$$
\begin{array}{lll}
  \integ{t}{\tau_j}&
\bigg(g_n(s,Y_{s}^{n})-g_m(s,Y_{s}^{m})\bigg)\psi'(Y_{s}^n
-Y_{s}^{m})1_{\{s\leq \tau_n\}}dA_s
\\ & \leq -\integ{t}{\tau_j}
g_m(s,Y_{s}^{m})\psi'(Y_{s}^n -Y_{s}^{m})1_{\{s\leq \tau_n\}}dA_s
\\ & \leq \integ{t}{\tau_j}
\psi'(Y_{s}^n -Y_{s}^{m})1_{\{s\leq \tau_n\}}dA_s
\\ & \leq 12j \integ{t}{\tau_j}e^{12j} (Y_s^n -Y_{s}^{m})dA_s,
\end{array}
$$
and
$$
\begin{array}{lll}
  \integ{t}{\tau_j}&
\bigg(f_n(s,Y_{s}^n,Z_{s}^n)
-f_m(s,Y_{s}^{m},Z_{s}^{m})\bigg)\psi'(Y_{s}^n -Y_{s}^{m})ds
\\ & \leq -\integ{t}{\tau_j}
f_m(s,Y_{s}^{m},Z_{s}^{m})\psi'(Y_{s}^n -Y_{s}^{m})1_{\{s\leq
\tau_n\}}ds
\\ & \leq \integ{t}{\tau_j}
(\eta_s+\frac{C_s}{2}|Z_{s}^{m}|^2)\psi'(Y_{s}^n
-Y_{s}^{m})1_{\{s\leq \tau_n\}}ds
\\ & \leq 12j \integ{t}{\tau_j}e^{12j} \eta_s(Y_s^n -Y_{s}^{m})ds
+\frac{j}{2} \integ{t}{\tau_j}|Z_{s}^{m}|^2\psi'(Y_{s}^n
-Y_{s}^{m})ds.
\end{array}
$$
It follows then that
\begin{equation}\label{eqq8}
\begin{array}{ll}
&\psi(Y_{t}^n -Y_{t}^{m}) +12 j \integ{t}{\tau_j}(Y_{s}^n -{L}_{s})
dK_{s}^{m+} +12 j\integ{t}{\tau_j}({U}_{s}-Y_{s}^m)
dK_{s}^{n-}\\&\leq \psi(Y_{\tau_j}^n
-Y_{\tau_j}^{m}) + 12 j \integ{t}{\tau_j}e^{12 j} 1_{\{\tau_n \leq s\leq \tau_m\}} dA_s\\
&+12j \integ{t}{\tau_j}e^{12j} (Y_s^n -Y_{s}^{m})dA_s+12 j
\integ{t}{\tau_j}e^{12j}{\eta}_s(Y_{s}^n -Y_{s}^{m})ds
-\integ{t}{\tau_j} dR_s^{n,m,j}\\ & +j
\integ{t}{\tau_j}\psi'(Y_{s}^n -Y_{s}^{m})|Z_{s}^n|^2 ds +j
\integ{t}{\tau_j}\psi'(Y_{s}^n -Y_{s}^{m})|Z_{s}^n -Z_{s}^m|^2 ds
\\ &
-\integ{t}{\tau_j}\psi'(Y_{s}^n -Y_{s}^{m})(Z_{s}^n-
Z_{s}^{m})dB_{s} -\dfrac12\integ{t}{\tau_j}\psi''(Y_{s}^n
-Y_{s}^{m})|Z_{s}^n- Z_{s}^{m}|^2 ds,
\end{array}
\end{equation} where $dR^{n, m,j}$ is a positive measure depending on $n, m$
and $j$. In order to understand the terms of $dR^{n, m,j}$, let us
give an example. Since $12j s\leq\psi'(s),\,\, \forall s\in [0,1 ]$,
the term
$$
-\integ{t}{\tau_j}\psi'(Y_{s}^n -{L}_{s}) dK_{s}^{m+},
$$
can be written as follows
$$
-\integ{t}{\tau_j}\psi'(Y_{s}^n -{L}_{s}) dK_{s}^{m+}
=-\integ{t}{\tau_j}\underbrace{\bigg(\psi'(Y_{s}^n -{L}_{s})
-12j(Y_{s}^n -{L}_{s})\bigg) dK_{s}^{m+}}_{\geq 0} -
\integ{t}{\tau_j}12j(Y_{s}^n -{L}_{s}) dK_{s}^{m+},
$$
the positive term in the above equation we put in $dR^{n, m,j}$, the
same holds for all estimates used in Equation (\ref{eqq8}). Now
$$ \begin{array}{ll}
& 2j \integ{t}{\tau_j}\psi'(Y_{s}^n -Y_{s}^{m})|Z_{s}^n -Z_{s}^m|^2
ds -\dfrac12\integ{t}{\tau_j}\psi''(Y_{s}^n -Y_{s}^{m})|Z_{s}^n-
Z_{s}^{m}|^2 ds \\ &= \integ{t}{\tau_j} (2j e^{12 j (Y_{s}^n
-Y_{s}^{m})} -2j -6je^{12 j (Y_{s}^n -Y_{s}^{m})})|Z_{s}^n-
Z_{s}^{m}|^2 ds,
\\ & = \integ{t}{\tau_j} -4j (e^{12 j (Y_{s}^n -Y_{s}^{m})} -1)|Z_{s}^n- Z_{s}^{m}|^2 ds -
6j \integ{t}{\tau_j}|Z_{s}^n- Z_{s}^{m}|^2 ds
\\ & = -4j \integ{t}{\tau_j}\psi'(Y_{s}^n
-Y_{s}^{m})|Z_{s}^n -Z_{s}^m|^2 ds -6j\integ{t}{\tau_j}|Z_{s}^n-
Z_{s}^{m}|^2 ds.
\end{array}
$$
Coming back to Equation (\ref{eqq8}) we obtain
\begin{equation}\label{eq33}
\begin{array}{ll}
&\psi(Y_{t}^n -Y_{t}^{m}) +12 j \integ{t}{\tau_j}(Y_{s}^n -{L}_{s})
dK_{s}^{m+} +12 j\integ{t}{\tau_j}({U}_{s}-Y_{s}^m) dK_{s}^{n-}
\\&+ 4j \integ{t}{\tau_j}\psi'(Y_{s}^n
-Y_{s}^{m})|Z_{s}^n -Z_{s}^m|^2 ds +6j\integ{t}{\tau_j}|Z_{s}^n-
Z_{s}^{m}|^2 ds +\integ{t}{\tau_j} dR_s^{n,m,j}
\\ & \leq \psi(Y_{\tau_j}^n -Y_{\tau_j}^{m}) + 12 j \integ{t}{\tau_j}e^{12
j} ((Y_s^n -Y_{s}^{m}) +1_{\{\tau_n \leq s\leq \tau_m\}}) dA_s \\ &+
12 j \integ{t}{\tau_j}e^{12j}{\eta}_s(Y_{s}^n -Y_{s}^{m})ds +2j
\integ{t}{\tau_j}\psi'(Y_{s}^n -Y_{s}^{m})|Z_{s}^n|^2 ds
\\ & -\integ{t}{\tau_j}\psi'(Y_{s}^n -Y_{s}^{m})(Z_{s}^n-
Z_{s}^{m})dB_{s}.
\end{array}
\end{equation} By taking $n=0$ in Equation (\ref{eq33}), inequality (\ref{eqq9}) and the fact that $Y^n -Y^m\leq 1$, we get
for all $j\in\N$
 $$
\begin{array}{ll}
& \E\integ{0}{\tau_j}|Z_{s}^0- Z_{s}^{m}|^2 ds \leq c_{j},
\end{array}
$$ where $c_j$ is a positive constant depending only on $j$.\\ In
force of inequality (\ref{eqq9}) we obtain, for all $j, m \in\N$
$$
\begin{array}{ll}
& \E\integ{0}{\tau_j}|Z_{s}^{m}|^2 ds \leq C_{j},
\end{array}
$$
where $C_j$ is a positive constant depending only on $j$.\\
Now, there exist a subsequence $\bigg(m_k^j\bigg)_{k}$ of $m$ and a
process $\widehat{Z}^j\in L^2(\Omega, L^2([0,\tau_j]; \R^d))$ such
that
\\
$Z_s^{m_k^j}1_{\{s\leq \tau_j\}}$ converges weakly in $L^2(\Omega,
L^2([0,\tau_j]; \R^d))$ to the process $\widehat{Z}_s^j1_{\{s\leq
\tau_j\}}$ as $k$ goes to infinity\\
 and
\\ $\sqrt{\psi'(Y_s^n -Y_s^{m_k^j})}(Z_s^n - Z_s^{m_k^j} )1_{\{s\leq
\tau_j\}}$ converges weakly in $L^2(\Omega, L^2([0,\tau_j]; \R^d))$
to $\sqrt{\psi'(Y_s^n -Y_s)}(Z_s^n - \widehat{Z}_s^{j} )1_{\{s\leq
\tau_j\}}$ as $k$ goes to infinity, since $\sqrt{\psi'(Y_s^n
-Y_s^{m_k^j})}$ converges strongly to $\sqrt{\psi'(Y_s^n -Y_s)}$.\\
Now, since
$$
\begin{array}{ll}
& 4j \E\integ{t}{\tau_j}\psi'(Y_{s}^n -Y_{s})|Z_{s}^n
-\widehat{Z}_{s}^j|^2 ds +6j\E\integ{t}{\tau_j}|Z_{s}^n-
\widehat{Z}_{s}^{j}|^2 ds\\&  \leq
\displaystyle{\liminf_{k}}\bigg[4j \E\integ{t}{\tau_j}\psi'(Y_{s}^n
-Y_{s}^{m_k^j})|Z_{s}^n -{Z}_{s}^{m_k^j}|^2 ds
+6j\E\integ{t}{\tau_j}|Z_{s}^n- {Z}_{s}^{m_k^j}|^2 ds\bigg],
\end{array}
$$
by taking $m = m_k^j,\,k\geq n$, in Equation (\ref{eq33}) and
tending $k$ to infinity, we get
 \begin{equation}\label{eq34}
\begin{array}{ll}
&\E\psi(Y_{t}^n -Y_{t})
  +4j \E\integ{t}{\tau_j}\psi'(Y_{s}^n
-Y_{s})|Z_{s}^n -\widehat{Z}_{s}^j|^2 ds
+6j\E\integ{t}{\tau_j}|Z_{s}^n- \widehat{Z}_{s}^{j}|^2 ds
\\ & \leq \E\bigg(\psi(Y_{\tau_j}^n -Y_{\tau_j})\bigg) + 12je^{12
j} \E\integ{t}{\tau_j} ((Y_s^n -Y_{s}) +1_{\{ s > \tau_n\}}) dA_s\\
&+12 j \E\integ{t}{\tau_j}e^{12j}{\eta}_s(Y_{s}^n -Y_{s})ds
 +4j \E\integ{t}{\tau_j}\psi'(Y_{s}^n
-Y_{s})|Z_{s}^n -\widehat{Z}_{s}^j|^2 ds \\ &+ 4j
\E\integ{t}{\tau_j}\psi'(Y_{s}^n -Y_{s})|\widehat{Z}_{s}^j|^2 ds.
\end{array}
\end{equation}
Hence $$\displaystyle\lim_{n\rightarrow +\infty}
\E\integ{0}{\tau_j}|Z_{s}^n- \widehat{Z}_{s}^{j}|^2 ds = 0.$$ 
By the uniqueness of the limit we obtain that $$
\widehat{Z}_s^j(\omega) 1_{\{0\leq s\leq \tau_j(\omega)\}} =
\widehat{Z}_s^{j+1}(\omega)1_{\{0\leq s\leq
\tau_j(\omega)\}},\,\,P(d\omega) ds-a.e.$$ For $s\in [0,T]$, let us
set $Z_s(\omega)= \displaystyle\lim_{j} \widehat{Z}_s^j 1_{\{ s \leq
\tau_j\}} = \widehat{Z}_s^{j(\omega)}(\omega),$ where $j(\omega)$ is
such that $\tau_{j(\omega)}(\omega) = T$. Then, for all $j\in\N$
$\E\integ{0}{\tau_j}\mid Z_s\mid^2 ds < +\infty$. Hence
$\integ{0}{T}\mid Z_s\mid^2 ds < +\infty, P-a.s. $ Moreover for all
$j\in\N$, we have
\begin{equation}\label{eqq10}
\displaystyle \lim_{n} \E\integ{0}{\tau_j}|Z_{s}^n -Z_{s}|^2 ds = 0.
\end{equation}
Proposition \ref{th1} proved.\eop
\subsection{Main Result} Now we are ready to give the main result of
this paper.
\begin{theorem}\label{the14}
Assume that assumptions $(\bf{H.0})$--$(\bf{H.4})$ hold. Then the
process $(Y, Z, K^+ , K^-)$, defined by (\ref{eqq5}) and Proposition
\ref{th1}, is the maximal solution for GRBSDE (\ref{eq221}).
\end{theorem}
\bop. Let us now prove that the process $(Y, Z, K^+ , K^-)$ is the
maximal solution for Equation (\ref{eq221}). To begin with, let us
show that the process $Y$ is continuous. From Equation (\ref{eq33})
and according to Bulkholder-Davis-Gundy inequality, there exists a
constant $C >0$ such that
\begin{equation}\label{eq35}
\begin{array}{ll}
 & \E\displaystyle\sup_{s\leq \tau_j}\psi(Y_{s}^n
-Y_{s})
\\ & \leq \E\bigg(\psi(Y_{\tau_j}^n -Y_{\tau_j})\bigg) + 12 je^{12
j} \E\integ{0}{\tau_j} ((Y_s^n -Y_{s}) +1_{\{ s > \tau_n\}}) dA_s \\
&+ 12j\E \integ{0}{\tau_j}e^{12j}{\eta}_s(Y_{s}^n -Y_{s})ds + 4j
\E\integ{0}{\tau_j}\psi'(Y_{s}^n -Y_{s})|Z_{s}|^2 ds
\\ &+C\E\bigg(\integ{0}{\tau_j}\mid\psi'(Y_{s}^n -Y_{s})\mid^2|Z_{s}^n
-Z_{s}|^2 ds\bigg)^{\frac12}.
\end{array}
\end{equation}
Hence $$ \lim_{n} \E\sup_{s\leq \tau_j}\psi(Y_{s}^n -Y_{s})=0,$$ and
then $$ \lim_{n} \E\sup_{s\leq \tau_j}\mid Y_{s}^n -Y_{s}\mid =0.$$
It
follows that $Y$ is continuous, since $P[\displaystyle{\cup_{j=1}^{\infty}(\tau_j= T)]} =1$.\\
Now, in view of (\ref{eqq10}) there exists a subsequence
$\bigg(n_k^j\bigg)_k$ of $n$ such that :
\begin{enumerate}
\item
$ \E\integ{0}{\tau_j}|Z_{s}^{n_k^j} -Z_{s}|^2 ds\leq
\frac{1}{2^k}\quad \mbox{and} \quad
\E\integ{0}{\tau_j}\sum_{k=0}^{+\infty}|Z_{s}^{n_k^j} -Z_{s}|^2
ds\leq 2, $
\item $Z_s^{n_k^j}(\omega)\longrightarrow Z_s(\omega),\, a.e.\,\,
(s,\omega)\in[0,\tau_j]\times \Omega,\quad\mbox{and}\quad \mid
Z_s^{n_k^j}(\omega)\mid \leq h^j_s,\,\, a.e.\,\,
(s,\omega)\in[0,\tau_j]\times \Omega,$ where $h_s^j = 1_{\{ s \leq
\tau_j\}}
\bigg(2|Z_s|^2+2\displaystyle\sum_{k=0}^{+\infty}|Z_s^{n_k^j}-Z_s|^2\bigg)^{\frac12}$.
\end{enumerate}
 Hence, in view of Lemma \ref{lem0} we obtain
$$
\begin{array}{ll}
& \E\integ{0}{\tau_j}| f_{n_k^j}(s, Y_s^{n_k^j}, Z_s^{n_k^j}) - f(s,
Y_s^{n_k^j}, Z_s^{n_k^j})| ds \\ &= \E\integ{0}{\tau_j} f_{n_k^j}(s,
Y_s^{n_k^j}, Z_s^{n_k^j}) - f(s, Y_s^{n_k^j}, Z_s^{n_k^j}) ds  \\ &
= \E\integ{0}{\tau_j}\bigg(f_{n_k^j}(s, Y_s^{n_k^j}, Z_s^{n_k^j}) -
f(s, Y_s^{n_k^j}, Z_s^{n_k^j})\bigg)
1_{\{ \mid Z_s^{n_k^j} -Z_s\mid \leq 1\}}ds\\
&\quad +\E\integ{0}{\tau_j}\bigg(f_{n_k^j}(s, Y_s^{n_k^j},
Z_s^{n_k^j}) - f(s, Y_s^{n_k^j}, Z_s^{n_k^j})\bigg) 1_{\{ \mid
Z_s^{n_k^j} -Z_s\mid \geq 1\}}ds
\\ &\leq  \E\integ{0}{\tau_j}\sup_{(y,z)\in [0,1]\times B(Z_s,1)}\bigg(f_{n_k^j}(s, y, z) -
f(s, y,z)\bigg) ds -\E\integ{0}{\tau_j} f(s, Y_s^{n_k^j},
Z_s^{n_k^j}) 1_{\{ \mid Z_s^{n_k^j} -Z_s\mid \geq 1\}}ds
\\ &\leq
\E\integ{0}{\tau_j}\sup_{(y,z)\in [0,1]\times
B(Z_s,1)}\bigg(f_{n_k^j}(s, y, z) - f(s, y,z)\bigg) ds
\\ &
\quad +\E\integ{0}{\tau_j}\bigg({\eta}_s +j \mid Z_s^{n_k^j}
-Z_s\mid^2 + j \mid Z_s\mid^2 \bigg)(\mid Z_s^{n_k^j} -Z_s\mid\wedge
1)ds,
\end{array}
$$
where $B(Z,1)$ is the closed ball of center $Z$ and radius $1$.\\
By taking account of Lemma \ref{lem0} one can see that
$$\displaystyle{\sup_{(y,z)\in [0,1]\times
B(Z_s,1)}\bigg(f_{n_k^j}(s, y, z) - f(s, y,z)\bigg)},$$ converges
pointwise. But, on $[0, \tau_j]$, we have also that
$$
\sup_{(y,z)\in [0,1]\times B(Z_s,1)}\bigg(f_{n_k^j}(s, y, z) - f(s,
y,z)\bigg) \leq \eta_s +\frac{j}{2}(|Z_s| +1)^2.
$$
Henceforth, by using Lebesgue's dominated convergence theorem, we
get
$$ \lim_k\E\integ{0}{\tau_j}\sup_{(y,z)\in [0,1]\times
B(Z_s,1)}\bigg(f_{n_k^j}(s, y, z) - f(s, y,z)\bigg) ds = 0, $$ and
$$\lim_k\E\integ{0}{\tau_j}\bigg({\eta}_s +j \mid Z_s^{n_k^j} -Z_s\mid^2 +
j \mid Z_s\mid^2 \bigg)(\mid Z_s^{n_k^j} -Z_s\mid\wedge 1)ds = 0.
$$
Therefore
$$ \lim_k \E\integ{0}{\tau_j}\bigg(f_{n_k^j}(s, Y_s^{n_k^j},
Z_s^{n_k^j}) - f(s, Y_s^{n_k^j}, Z_s^{n_k^j})\bigg) ds =0.$$ It
follows then from Lebesgue's dominated convergence theorem that for
all $j\in\N$
 $$  \lim_k \E\integ{0}{\tau_j}\mid f(s, Y_s^{n_k^j},
Z_s^{n_k^j}) - f(s, Y_s, Z_s)\mid ds =0.$$ Hence for all $j\in\N$ $$
\lim_k \E\integ{0}{\tau_j}\mid f_{n_k^j}(s, Y_s^{n_k^j},
Z_s^{n_k^j}) - f(s, Y_s, Z_s)\mid ds =0.$$ Since the above limit
doesn't depend on the choice of the subsequence $(n_k^j)_{k}$ we
have for all $j\in\N$
$$  \lim_n
\E\integ{0}{\tau_j}\mid f_{n}(s, Y_s^{n}, Z_s^{n}) - f(s, Y_s,
Z_s)\mid ds =0.$$ It not difficult also to prove that for all
$j\in\N$
$$  \lim_n
\E\integ{0}{\tau_j}\mid g_{n}(s, Y_s^{n}) - g(s, Y_s)\mid dA_s =0.
$$
From Equation (\ref{eqq})$(i)$ we obtain, $\forall j$,\, $
 \sup_{n}\E K_{\tau_j}^{n+} <+\infty$. It then follows from Fatou's lemma that for any $j\in\N$,\, $ \E
K_{\tau_j}^{+}<+\infty$. Henceforth $K_{T}^{+}< +\infty, P$-a.s.
\\ Then we obtain, $P-$a.s.
$$
\begin{array}{ll}
Y_{t}= & Y_{\tau_j} + \integ{t}{\tau_j}f(s,Y_{s},Z_{s})ds+
\integ{t}{\tau_j} g(s,Y_s)dA_s + \integ{t}{\tau_j}dR_s
\\ & +\integ{t}{\tau_j}dK_{s}^{+} -\integ{t}{\tau_j}dK_{s}^{-}
-\integ{t}{\tau_j}Z_{s}dB_{s},
\end{array}
$$
Since $\tau_j$ is a stationary stopping time we get $P-$a.s.
$$
\begin{array}{ll}
Y_{t}= & \xi + \integ{t}{T}f(s,Y_{s},Z_{s})ds+ \integ{t}{T}
g(s,Y_s)dA_s + \integ{t}{T}dR_s
\\ & +\integ{t}{T}dK_{s}^{+} -\integ{t}{T}dK_{s}^{-}
-\integ{t}{T}Z_{s}dB_{s},
\end{array}
$$

Now let us prove the minimality conditions. We have
$$
\integ{0}{T}( {U}_{t}-Y_{t}^{n}) dK_{t}^{n-}=0.
$$
Hence, since $dK^{-} = \displaystyle\inf_{n}dK^{n-}$, we get
$$
\integ{0}{T}( {U}_{t}-Y_{t}^{n}) dK_{t}^{-}=0.
$$
 It follows then from Fatou's lemma that
$$
\integ{0}{T}( {U}_{t}-Y_{t}) dK_{t}^{-}  = 0.
$$
On the other hand
$$
\integ{0}{T}(Y_{t}^{n} -L_t) dK_{t}^{n+}  = 0.
$$
Hence, since $Y = \displaystyle\inf_{n}Y^{n}$, we obtain
$$
\integ{0}{T}(Y_{t} -L_t) dK_{t}^{n+} = 0.
$$
Applying Fatou's lemma we obtain
$$
\integ{0}{T}(Y_{t} -L_t) dK_{t}^{+} = 0.
$$
Now, since $dK^{+} = \displaystyle\sup_{n}dK^{n+}$,\,\, $dK^{-} =
\displaystyle\inf_{n}dK^{n-}$ and the measures $dK^{n+}$ and
$dK^{n-}$ are singular, it follows that $dK^{+}$ and $dK^{-}$ are
singular.
\\ Now it is not difficult to see that the process $(Y, Z, K^+ ,
K^-)$ satisfies Equation (\ref{eq221}). It remains to prove $(Y, Z,
K^+ , K^-)$ is maximal. Let $(Y', Z', K^{'+} , K^{'-})$ be another
solution to Equation (\ref{eq221}). By comparison theorem we have
that $Y' \leq Y^n$ and then $Y'\leq Y$. The proof of Theorem
\ref{th1} is then finished.\eop \subsection{Existence of maximal
solution for Equation (\ref{eq0})} Given the result of Theorem
\ref{the14} and taking advantage of Proposition \ref{pro3}, the
following theorem follows directly by a logarithmic change.
\begin{theorem}\label{the2}
Let assumptions $(\bf{A.1})$--$(\bf{A.3})$ hold true. Then there
exists a maximal solution for Equation (\ref{eq0}).
\end{theorem}

\bop.  Let $(\overline{Y}_t ,\overline{{Z}}_t ,\overline{{K}}^+_t
,\overline{{K}}^-_t )_{t\leq T}$ be the maximal solution of Equation
(\ref{eq221}) then, for any $t\leq T$, we have
\begin{equation}
\label{eq222} \left\{
\begin{array}{ll}
(i) & 
 \overline{Y}_{t}=\overline{\xi} +\integ{t}{T}\overline{f}(s,\overline{Y}_{s},\overline{Z}_{s})ds +
 \integ{t}{T}g(s,\overline{Y}_{s})dA_s + \integ{t}{T} dR_s\\
&\qquad\quad+\integ{t}{T}d\overline{K}_{s}^+
-\integ{t}{T}d\overline{K}_s^-
-\integ{t}{T}\overline{Z}_{s}dB_{s}\,, t\leq T,
\\ (ii)&
\forall  t\leq T,\,\, \overline{L}_t \leq \overline{Y}_{t}\leq
\overline{U}_{t},\\ (iii)  & \integ{0}{T}(
\overline{Y}_{t}-\overline{L}_{t}) d\overline{K}_{t}^+=
\integ{0}{T}( \overline{U}_{t}-\overline{Y}_{t})
d\overline{K}_{t}^-=0,\,\, \mbox{a.s.} \\ (iv)& \overline{Y}\in
{\cal C} \quad \overline{K}^+, \overline{K}^-\in {\cal K} \quad
\overline{Z}\in {\cal L}^{2,d},
 \\ (v)&
d\overline{K}^+\perp d\overline{K}^-.
\end{array}
\right.
\end{equation}
Now, for all $t\leq T$, let us set $$ Y_t =
\dfrac{\ln(\overline{Y}_t)}{m_t} +m_t, \quad  Z_t =
\dfrac{\overline{Z}_t}{m_s \overline{Y}_s }, \quad dK_t^{\pm} =
\dfrac{d\overline{K}_t^{\pm}}{m_t \overline{Y}_t }.$$ By using
It\^{o}'s formula to $\dfrac{\ln(\overline{Y}_t)}{m_t} +m_t$, we
have $$
 Y_{t}=\xi
+\integ{t}{T}f(s,Y_{s},Z_{s})ds +\integ{t}{T}dK_{s}^+
-\integ{t}{T}dK_{s}^- -\integ{t}{T}Z_{s}dB_{s}\,, t\leq T.
$$
Therefore it is not difficult to prove that $(Y, Z, K^+, K^-)$ is a
maximal solution for the GBSDE with two reflecting barriers
(\ref{eq0}). This completes the proof. \eop

\section{\Large{\bf{Appendix}}}
\appendix
\section{Comparison theorem }
The comparison theorem for real-valued BSDEs turns out to be one of
the classic results of the theory of BSDE. It allows to compare the
solutions of two real-valued BSDEs whenever we can compare the
terminal conditions and the generators. This section is devoted to
present a comparison theorem for the following GBSDE with generator
$hdA_s$:
\begin{equation}
\label{eq22111} \left\{
\begin{array}{ll}
(i) &
 {Y}_{t}={\xi} +\integ{t}{T}h(s,{Y}_{s},{Z}_{s})dA_s
+\integ{t}{T}d{K}_{s}^+ -\integ{t}{T}d{K}_s^-
-\integ{t}{T}{Z}_{s}dB_{s}\qquad t\leq T,
\\ (ii)&
\forall  t\leq T,\,\, {L}_t \leq {Y}_{t}\leq {U}_{t},\\ (iii)
&\integ{0}{T}({Y}_{t}-{L}_{t}) d{K}_{t}^+= \integ{0}{T}(
{U}_{t}-{Y}_{t}) d{K}_{t}^-=0,\,\, \mbox{a.s.}
\\ (iv)& {Y}\in {\cal C} \quad {K}^+,
{K}^-\in {\cal K} \quad {Z}\in {\cal L}^{2,d},
 \\ (v)&
d{K}^+\perp d{K}^-.
\end{array}
\right. \end{equation} Let $(Y^i, Z^i, K^{i+}, K^{i+})$ $(i=1,2)$ be
two solutions (if they exist) of Equation (\ref{eq22111}) associated
respectively with $(\xi^1, h^1, A^1, L^1, U^1)$ and $(\xi^2, h^2,
A^2, L^2, U^2)$, such
that, for $(i=1,2)$ the following assumptions are satisfied :\\
\noindent$(\bf{D.1})$  $L^i: [0,T]\times \Omega\longrightarrow
\R\cup\{-\infty\}$ and  $U^i: [0,T]\times \Omega\longrightarrow
\R\cup\{+\infty\}$ are two continuous barriers processes satisfying
$$
L^1_t\leq L^2_t\quad\mbox{and}\quad U^1_t\leq U^2_t, \forall t\in
[0,T].
$$
 \noindent$(\bf{D.2})$ $\xi^1 \leq \xi^2$.\\
 \noindent$(\bf{D.3})$ $A^i\in \cal K$ and for all $(s, \omega),\quad$ $ h^1(s,{Y}_{s}^1,{Z}_{s}^1)dA_s^1\leq
h^2(s,{Y}_{s}^1,{Z}_{s}^1)dA_s^2 $.\\
 \noindent$(\bf{D.4})$ There exist two processes $\alpha \in \cal K$ and $b\in {\cal L}^{2,1}$ such that :
 $$ (h^2(s,{Y}_{s}^1,{Z}_{s}^1)-h^2(s,{Y}_{s}^2,{Z}_{s}^2))dA_s^2\leq
\mid Y_s^1 -Y_s^2\mid d\alpha_s + \mid Z_s^1 -Z_s^2\mid b_s ds.$$
\\
Set $\Gamma_t = e^{\widetilde{\alpha}_t + \int_0^t\widetilde{b}_s
dB_s -\frac12 \int_0^t\mid \widetilde{b}_s\mid^2 ds}$, where the
processes $\widetilde{\alpha}$ and $\widetilde{b}$ are defined as
follows :
$$
\widetilde{\alpha}_s = \dint_0^s \dfrac{Y^1_t -Y^2_t}{\mid Y^1_t
-Y^2_t\mid }1_{\{Y^1_t \neq Y^2_t\}}d\alpha_t\,\,\,
\mbox{and}\,\,\,\widetilde{b}_t = b_t \dfrac{Z^1_t -Z^2_t}{\parallel
Z^1_t -Z^2_t\parallel}1_{\{Z^1_t \neq Z^2_t\}}.
$$
We have the following. For the proof see \cite{EH}.
\begin{theorem} (Comparison theorem) \label{th13} Assume that assumptions $(\bf{D.1})-(\bf{D.4})$ hold.
\begin{itemize}
\item[i)] If\, $\displaystyle\liminf_{r\rightarrow +\infty} r\,
P\bigg[\sup_{0\leq s\leq T}\Gamma_s (Y^1_s -Y^2_s)^+>r\bigg] =0,$ we
have $Y^{1}_t\leq Y^{2}_ t$, for $t\in [0,T]$, a.s.
\item[ii)] If\, $Y^{1}_t\leq
Y^{2}_ t$, for $t\in [0,T]$, a.s., then :
$$
1_{\{U^1_t = U^2_t\}}dK^{1-}_t \leq dK^{2-}_t\,\,\, \mbox{and}\,\,\,
1_{\{L^1_t = L^2_t\}}dK^{2+}_t \leq dK^{1+}_t.
$$
\end{itemize}
\end{theorem}
We give the following remarks.
\begin{remark} It should be noted that :
\begin{enumerate}
\item If $\E\displaystyle\sup_{t\leq T} \Gamma_t (Y^1_t
-Y^2_t)^+<+\infty$ then $\displaystyle\lim_{r\rightarrow +\infty}
r\, P\bigg[\sup_{0\leq s\leq T}\Gamma_s (Y^1_s -Y^2_s)^+>r\bigg]
=0.$
\item If there exist constants $C>0$ and $p>1$  such that $\alpha_T +
\dint_0^T \mid b_s\mid^2 ds \leq C$ and $\E\displaystyle\sup_{t\leq
T} ((Y^1_t -Y^2_t)^{+})^p <+\infty,$ then we obtain
$\E\displaystyle\sup_{t\leq T} \Gamma_t (Y^1_t -Y^2_t)^+<+\infty$.
\item If  $U^{i} \equiv +\infty$, then $dK^{i-} \equiv 0,$ for $i=1, 2$.\\
\item If  $L^{i} \equiv -\infty$, then $dK^{i+} \equiv 0,$ for
$i=1, 2$.
\end{enumerate}
\end{remark}

\section{Existence and uniqueness of solutions for GRBSDE under strong assumptions on the coefficients}
 In this section we shall present the existence of solutions to
GRBSDE (\ref{eq0}) under strong assumptions on the coefficients.
\\ \\  We assume the following assumptions :
\\ $\bf{(C.1)}$ $(i)$ $f$ is uniformly Lipschitz with respect to
$(y,z)$, \it{i.e.}, there exists a constant $0<C_1<\infty$ such that
for any
 $y$, $y^{\prime }$,
$z$, $z^{\prime }\in \R$,
$$
|f( \omega ,t,y,z) -f(\omega ,t,y',z')| \leq C_1( | y-y'| +|z-z'|).
$$
$(ii)$ There exists a constant $C_2>0$ such that for all $y$,
$y^{\prime }$,
$z$, $z^{\prime }\in \R$, $-C_2\leq f(\omega ,t,y,z)\leq 0$.\\
$\bf{(C.2})$  $(i)$ $g$ is uniformly Lipschitz with respect to $y$,
\it{i.e.}, there exists a constant $0<C_3<\infty$ such that for any
 $y$, $y^{\prime }\in \R$,
$$
|g( \omega ,t,y) -g(\omega ,t,y')| \leq C_3 | y-y'|.
$$
$(ii)$ For all $t\in [0, T]$,\,\, $y\in \R$, $-1\leq g(t,y)\leq 0$.\\
$\bf{(C.3})$ For all $t\in [0,T],\,\,$ $0 <L_t\leq U_t< 1$.
\\
$\bf{(C.4})$ There exist constants $C_4, C_5>0$ such that $A_T \leq C_4$ and $\mid R\mid_T \leq C_5$.\\
\begin{theorem}\label{th22} Let assumptions $\bf{(C.1})-\bf{(C.4})$ and $\bf{(H.4})$ hold. Then there exists a
unique solution for GRBSDE (\ref{eq0}). Moreover,
$$\E\dint_0^{T}\mid Z_s \mid^2ds+ \E(K^{\pm}_ {T})^2 < +\infty.$$
\end{theorem}
\bop. Uniqueness follows directly from comparison theorem. For the
existence proof see \cite{EH}.\\ \\
%
%
%
\noindent\bf{Acknowledgment. } The authors would like to thank the
anonymous referee for his comments that help improve the paper.

\end{document}